\documentclass{article}

\usepackage{graphicx}
\usepackage{amsmath,amsthm}
\usepackage{amssymb,mathtools,mathabx,bm}
\usepackage[round]{natbib}
\usepackage{tikz,tikz-cd}
\usepackage[affil-it]{authblk}
\usepackage[letterpaper,twoside,vscale=.75,hscale=.70,nomarginpar,hmarginratio=1:1]{geometry}
\usepackage{setspace}
\usepackage[pageanchor=true,plainpages=false, pdfpagelabels, bookmarks,bookmarksnumbered]{hyperref}

\newtheorem{theorem}{Theorem}
\numberwithin{theorem}{section}
\newtheorem{proposition}[theorem]{Proposition}
\newtheorem{lemma}[theorem]{Lemma}

\newtheorem{corollary}[theorem]{Corollary}

\newtheorem{defn}[theorem]{Definition}
\newtheorem{thm}[theorem]{Theorem}
\newtheorem{prop}[theorem]{Proposition}
\newtheorem{lem}[theorem]{Lemma}

\newtheorem{rem}[theorem]{Remark}

\newtheorem{example}[theorem]{Example}

\newcommand*{\cat}[1]{{\mathbf{#1}}} 
\newcommand{\D}{\cat{D}}
\newcommand{\C}{\cat{C}}
\newcommand{\E}{\mathcal{E}}
\newcommand{\F}{\mathcal{F}}
\newcommand*{\Cat}{\cat{Cat}}
\newcommand*{\Set}{\cat{Set}}

\newcommand*{\Pos}{\cat{Pos}} 
\newcommand*{\Topos}{\cat{Topos}}
\newcommand*{\PshC}{\cat{Set}^{\C\op}}
\newcommand*{\PshD}{\cat{Set}^{\D\op}}

\newcommand*{\IndCat}{\cat{IndCat}}

\newcommand*{\Comon}{\cat{Comon}}

\newcommand*{\N}{{\mathbb{N}}}

\newcommand*{\monoto}{\rightarrowtail}

\newcommand*{\inc}{\hookrightarrow}

\newcommand*{\inv}{^{-1}}
\newcommand*{\op}{^{\mathrm{op}}}

\DeclareMathOperator{\Hom}{Hom}

\DeclareMathOperator{\id}{id}
\DeclareMathOperator{\Id}{id}

\DeclareMathOperator{\sub}{Sub}
\DeclareMathOperator{\Sh}{Sh}
\DeclareMathOperator{\St}{St}

\DeclareMathOperator{\y}{y}

\DeclareMathOperator{\Coalg}{Coalg}

\DeclareMathOperator{\Inc}{Inc}

\DeclareMathOperator{\RT}{RT}
\DeclareMathOperator{\PAss}{PAss}

\newcommand{\ve}{\varepsilon}
\newcommand{\ar}{\arrow}
\newcommand{\f}{\Box}

\newcommand{\mcE}{\mathcal{E}}

\newcommand{\ie}{\emph{i.e.}}
\newcommand{\eg}{\emph{e.g.}}

\newcommand{\cf}{\emph{cf.}}
\newcommand{\etc}{\emph{etc.}}

\newcommand{\ifoif}{if and only if}
\newcommand*{\ko}{k}

\DeclarePairedDelimiter\abs{\lvert}{\rvert}

\title{Stratified Toposes}

  \author{Colin Zwanziger}
  \affil{Institute of Philosophy, Czech Academy of Sciences}

\begin{document}
\onehalfspacing

\maketitle
\begin{abstract}
    We introduce \textbf{stratified toposes}, which are toposes that are stratified by a suitable hierarchy of universes. The term `stratified topos' recalls the notion of stratified pseudotopos of \citet{moerdijk2002type}. However, the details of our proposal are closer to one of \citet{streicher2005universes}, with the foundational contribution being a strengthening of Streicher's axioms. As such, stratified toposes model the calculus of constructions. 

   Key results about toposes can be refined to yield results about stratified toposes. As proof of concept, we construct what we call the \textbf{stratified topos of coalgebras} for a stratified Cartesian comonad on a stratified topos. 
    This construction refines that of the topos of coalgebras for a Cartesian comonad on a topos. Coalgebra constructions in related settings were given by \citet{warren2007coalgebras} and \citet{zwanziger2023natural}.

    This coalgebra construction exposes conditions under which generic monomorphisms and dense subcategories can be lifted to coalgebras, which are extracted in the appendix.
\end{abstract}

\section{Introduction}

The notion of (elementary) topos \citep{lawvere1970lectures} abstracts to the level of categorical algebra several aspects of the category of sets. However, it is natural to assume the existence in the category of sets of large cardinals, which is not reflected in the topos axioms. 

In remedy of this, various notions of universe in a topos have been introduced. \citet{benabou1973problemes}’s axioms are already quite close to the present approach. The axioms of \citet{streicher2005universes} are stronger than B\'enabou’s, and draw on previous work on the semantics of the calculus of constructions \citep[\cf][]{hyland1989theory,streicher1991semantics}. The present approach involves 
a novel strengthening of Streicher’s axioms. Compatible principles isolated in work on algebraic set theory since \citet{joyal1995algebraic} could potentially be added in the future.

We will emphasize that a topos which is stratified by a hierarchy of universes, which we call a \textbf{stratified topos}, is a natural categorical gadget. Whereas, in an ordinary topos, monomorphisms are `represented' by the subobject classifier, in a stratified topos, \textit{all morphisms} are, moreover, `represented' by some universe. As long as the notion of universe used in a given definition of stratified topos is reasonable, much of the category theory of toposes will be replicable at the level of stratified toposes. This viewpoint is evident in \citet{moerdijk2002type}, though their relatively baroque axioms for what they call a `stratified pseudotopos' may have inhibited broad adoption of their insight.

Of course, the present approach to universes is reasonable enough to yield such a well-behaved notion of stratified topos. In hewing close to \citet{streicher2005universes}, our choice of axioms is motivated by the view that stratified toposes should provide a categorical model for the calculus of constructions, in much the same sense that toposes provide a categorical model for constructive logic. At a minimum, any Grothendieck topos will contain a stratified topos, provided enough Grothendieck universes are assumed in the underlying set theory. We will not answer the admittedly important question of whether  realizability toposes yield stratified toposes in our sense, since this depends on an open problem.

To motivate the present approach to universes and stratified toposes, let us return to the set-theoretic case. When we assume, as in Tarski-Grothendieck set theory, that every set is contained in a Grothendieck universe, we can regard the class $\mbox{Set}$ of all sets as stratified by a hierarchy of Grothendieck universes
\[\mathrm{Set}_0 \subset \mathrm{Set}_1 \subset … \mathrm{Set}_n \subset ….\]
Moreover, we have 
\[\mathrm{Set}_0 \in \mathrm{Set}_1 \in … \mathrm{Set}_n \in ….\]
Grothendieck universes thus live a double life as subset and as element.

 In a categorical setting, it is natural to think of universes themselves as categories. Thus, in our notion of stratified topos, we will ask for a sequence of toposes and logical inclusions
    \[\mcE_0 \inc \mcE_1 \inc …\mcE_n \inc …,\]
such that, whenever $m<n$, $\mcE_n$ admits an internal category $\pmb{\mho}_m  \in \Cat(\mcE_n)$ that `represents’  $\mcE_m$, in the sense that 
\[\Hom_{\Cat(\mcE_n)}(I,\pmb{\mho}_m) \simeq (\mcE_n/I)_{<\E_m},\]
pseudonaturally in $I \in \mcE_n$, where $(\mcE_n/I)_{<\E_m}$ is the full subcategory of $\mcE_n/I$ on the maps with `$\mcE_m$-small fibers.’ (It is too much to ask that `$\mcE_m \in \Cat(\mcE_n)$' in general.) More intuitively, we then have, in particular, that \[\Hom_{\Cat(\mcE_n)}(I,\pmb{\mho}_m) \simeq \mcE_m/I,\]
pseudonaturally in $I \in \mcE_m$.

We will also ask for a novel density condition on the logical inclusion \[\mcE_m \xhookrightarrow{i} \mcE_n,\] which is used to extend the induced $\E_m$-indexed logical inclusion  
\[\mcE_m/(-) \inc i^* \mcE_n/(-)\]
to the expected $\E_n$-indexed logical inclusion
\[(\mcE_n/(-))_{<\E_m} \inc \mcE_n/(-).\]  This will be enough to ensure that our universes satisfy Streicher's axioms.\footnote{More informally, we note that any density condition  asserts a kind of similarity of $\mcE_n$ to $\mcE_m$, and is thus akin in motivation to a set-theoretic reflection principle. However, a reflection principle would assert a similarity of $\mcE_m$ to $\mcE_n$. An apt term for such a converse reflection principle would be \textbf{projection principle}. The described use case for our density axiom involves `projecting' properties like Cartesian closure from the $\mcE_m$-indexed category $\mcE_m/(-)$ to the $\mcE_n$-indexed category $(\mcE_n/(-))_{<\E_m}$.}

Towards showing that the tidy notion of stratified topos introduced is well-behaved, we will show that, like toposes, these stratified toposes are closed under the construction of coalgebras for suitable comonads. We introduce a suitable notion of \textbf{stratified Cartesian functor} as the analog of the notion of Cartesian functor in topos theory. This allows us to define a notion of \textbf{stratified Cartesian comonad}. We then construct the \textbf{stratified topos of coalgebras} for such a stratified Cartesian comonad. 
 This natural construction solves in the present setting a problem that was left open by \citet{moerdijk2002type} in the setting of stratified pseudotoposes, and solved in  \citet{streicher2005universes}'s setting in the dissertation of the author \citep{zwanziger2023natural}, and another setting by \citet{warren2007coalgebras}.

 Section \ref{sec:univ} introduces our notion of stratified topos, Section \ref{sec:univ'} constructs the stratified topos of coalgebras, and Section \ref{sec:conc} offers parting thoughts. Appendix \ref{sec:appendix} isolates some coalgebra constructions of general interest (lifting generic monos and dense subcategories to coalgebras).

\section{Universes and Stratified Toposes}\label{sec:univ}

In this section, we introduce our notion of \textbf{universe} and the resulting notion of \textbf{stratified topos}. Our notion of universe is based upon a novel notion of \textbf{dense logical subtopos}.

 Section \ref{sec:conv} recalls some categorical notation and fundamentals. Section \ref{sec:streicher} reviews the universes of \citet{streicher2005universes}. Section \ref{sec:unibasic} starts to develop our approach to universes, motivating dense logical subtoposes. Section \ref{sec:embden'} contains a necessary divagation on dense functors. Section \ref{sec:denslogsub} introduces dense logical subtoposes. Finally, Section \ref{sec:strattop} introduces our notion of stratified topos.

\subsection{Conventions}\label{sec:conv}

When $\C$ is a category, we write  $\abs{\C}$ for its class of objects and $\mathrm{C}_1$ for its class of arrows.

We will assume that a Cartesian category comes equipped with a choice of finite limits. When $\D$ is a Cartesian category and $\C$ 
is a full subcategory of $\D$, let us say that $\C$ is \textbf{Cartesian subcategory} when it is closed under the choice of finite limits in $\D$.

We denote by $\Cat$ the $2$-category of small categories and by $\cat{CAT}$ the (very large) $2$-category of categories. We denote by $\cat{Topos}$ the $2$-category of small toposes, Cartesian functors, and natural transformations. 

When $\C$ is a category and $\D$ is a $2$-category, we denote by $\D^{\C}$ the $2$-category of $2$-functors from $\C$ (regarded as a locally discrete $2$-category) to $\D$, $2$-natural transformations, and modifications.

By a \textbf{$\C$-indexed category} we mean a contravariant pseudofunctor from a category $\C$ (regarded as a $2$-category) to $\cat{CAT}$. By a \textbf{$\C$-indexed functor} (resp. \textbf{natural transformation}) we mean a pseudonatural transformation (resp. modification) between  $\C$-indexed categories (resp. functors). 

For the rest of the section we assume a small category $\C$. 

We say that a $\C$-indexed category $P$ is \textbf{small} if $P(I)$ is a small category, for all $I \in \C$. 
We denote by $\IndCat(\C)$ the $2$-category of small $\C$-indexed categories and $\C$-indexed functors and natural transformations.   When this $J$ is a Grothendieck coverage on $\C$, we denote by $\St(\C, J)$ the sub-$2$-category of $\IndCat(\C)$ consisting of the (small) stacks on $(\C, J)$ and $\C$-indexed functors and natural transformations.

We write $\Cat(\C)$ for the $2$-category of categories in $\C$. When $\cat{D}$ is a category in $\C$, we write $[\cat{D}] \in \Cat^{\C\op}$ for its \textbf{externalization}, given on $I \in \C$ by \[[\D](I) = \Hom_{\Cat(\C)}(I,\cat{D}),\]
where $I$ is regarded on the right as a discrete internal category. We say that $\cat{D}$  \textbf{internalizes} some functor $P : \C\op \to \Cat$ if we have $[\cat{D}] \cong \cat{P} \in \Cat^{\C\op}$. We recall that the operation 
\[\D \mapsto [\D]\]
extends to a $2$-functor 
\[[-] : \Cat(\C) \to \Cat^{\C\op},\] 
which is a $2$-embedding in the sense that it induces an isomorphism of categories 
\[\Hom_{\Cat(\C)}(\D,\cat{E}) \cong \Hom_{\Cat^{\C\op}}([\D],[\cat{E}]),\]
for any categories $\D$ and  $\cat{E}$ in $\C$ \citep[see][\S 1]{street2018categories}. To generalize our initial usage, when $D : \cat{J} \to \Cat(\C)$ is a $2$-diagram (\ie{} a $2$-functor), we  say that the composite $[-]\circ D$ is the \textbf{externalization} of $D$, and that such a $D : \cat{J} \to \Cat(\C)$ \textbf{internalizes} some $2$-diagram $D' : \cat{J} \to \Cat^{\C\op}$ if $[-]\circ D \cong D'$.\footnote{In this connection, the only nontrivial $2$-diagrams we will be interested in are those that encode comonads (appearing in Section \ref{sec:univ'}).}

\subsection{Streicher Universes}\label{sec:streicher}

Since it will be a point of reference for us, we start by introducing the topos-theoretic universes of \citet{streicher2005universes}. We include some motivation for the stronger axioms introduced later in  Section \ref{sec:univ}.

\begin{defn}[\citealt{taylor1999practical}]
When $\D$ is a category, a \textbf{class of display maps} $\mathrm{C}$ consists of a class of morphisms of $\D$, called \textbf{display maps}, 
such that    
\begin{itemize}
  \item $\mathrm{C}$ is closed under composition with isomorphisms;
  \item for every display map $p: E \to B$ and morphism $f : A \to B$, the pullback $f^*p : f^* E \to A$ exists and is a display map.
    \end{itemize}
\end{defn}

Given a class of display maps $\mathrm{C}$ in a category $\D$, 
we write $\mathbb{C}(A)$ for the full subcategory of $\D/A$ spanned by the display maps with codomain $A\in \D$. We write $\mathbb{C}$ for the $\D$-indexed category with value $\mathbb{C}(A)$ for each $A\in \D$ and pseudofunctorial action given by pullback. We say that a display map $p : \mho_\bullet \to \mho$ is \textbf{generic} if, whenever $f : A \to B$ is a display map, there exists a pullback square of form
\begin{equation*}
\begin{tikzcd}
A\ar[dd, "f"' ]\ar[ddrr, phantom, "\lrcorner", very near start]\ar[rr, ""]&& \mho_\bullet \ar[dd, "p"]\\
\\
B \ar[rr, ""'] &&\mho
\end{tikzcd}
\end{equation*}
in $\D$.

The following formulation is due to Mike Shulman \citep{shulman2009universe}:

\begin{defn}\label{def:imuniv}
     A class of display maps $\mathrm{E}$ in a topos $\F$ is said to be an \textbf{impredicative universe} (in the sense of \citet{streicher2005universes}) if 
     \begin{itemize}
         \item the composition of two display maps is again a display map;
        \item $\mathrm{E}$ admits a generic display map;
        \item our indexed subcategory $\mathbb{E}$ of $\F/(-)$ induced by $\mathrm{E}$ is an indexed  logical subtopos.
     \end{itemize}
\end{defn}

 When $\mathrm{E}$ is such an impredicative universe in a topos $\F$, we refer to display maps rather as ($\mathrm{E}$\textbf{-})\textbf{small maps}. We say that an object $A \in \F$ is \textbf{small} if the morphism $! : A \to 1$ is a small map.

 While Streicher's notion is already adequate for most purposes, it may reasonably be strengthened by further assumptions. For example, given such an impredicative universe $\mathrm{E}$ in a topos $\F$, let us say that a morphism $p : E \to B$ of $\F$ is a \textbf{locally small map}, 
 if, for any small $A \in \F$ and $f : A \to B$, the object $f^*E$ is small.

 Then, we have the following:
 \begin{prop}\label{prop:smallsmall}
     Any small map is a locally small map.
 \end{prop}

 \begin{proof}
     Given a small map $p : E \to B$, small object $A$, and morphism $f : A \to B$, the unique morphism $f^*E \to 1$ arises as the composite of small maps \[f^*E \xrightarrow{f^*p} A \to 1.\]
 \end{proof}

The converse of Proposition \ref{prop:smallsmall} is not derivable. When it does hold, $\mathbb{E} \inc \F/(-)$ can be recovered from  $\mathbb{E}(1) \monoto \F$, in line with the set-theoretic intuition  that $\mathbb{E}(1)$ should determine the notion of smallness in  $\F$. 

This converse could simply be assumed as an axiom. However, we will prefer to derive it from a stronger set of axioms, justifiable on independent grounds.

\subsection{Our Basic Approach to Universes}\label{sec:unibasic}

We now take the first steps towards our axiomatization of universes in toposes, also motivating the addition of our novel density axiom in Section \ref{sec:denslogsub}.

As discussed in the introduction, our approach will be based on the notion of logical inclusion of toposes 
\[\E \inc \F\]
(rather than the less fundamental notion of indexed logical inclusion of form \[\mathbb{E} \inc \F/(-)\]
that is traditionally considered, as in Section \ref{sec:streicher}). A universe in this context will, morally, be a category $\pmb{\mho}$ in $\F$ which suitably `represents' $\E$.

More precisely, when $\C$ is a Cartesian full subcategory of a Cartesian category $\D$, we say that a morphism $p : E \to B$ of $\D$ is a \textbf{locally} ($\C$\textbf{-})\textbf{small map}, if, for any $A \in \C$ and $f : A \to B$, we can form a pullback
\begin{equation*}
\begin{tikzcd}
D\ar[dd, ""' ]\ar[ddrr, phantom, "\lrcorner", very near start]\ar[rr, ""]&& E \ar[dd, "p"]\\
\\
A \ar[rr, "f"'] &&B,
\end{tikzcd}
\end{equation*}
 in which $D \in \C$. For convenience, we assume that when such a $p$ is locally small and we have $A \in \C$ and $f : A \to B$, the chosen pullback
\begin{equation*}
\begin{tikzcd}
f^* E\ar[dd, "f^*p"' ]\ar[ddrr, phantom, "\lrcorner", very near start]\ar[rr, ""]&& E \ar[dd, "p"]\\
\\
A \ar[rr, "f"'] &&B
\end{tikzcd}
\end{equation*}
has $f^* E \in \C$. 

 Locally small maps are closed under pullback, and thus form a class of display maps in $\D$. We denote by $(\D/I)_{<\C}$ the full subcategory of $\D/I$ spanned by the locally small maps with codomain $I \in \D$, and by $(\D/(-))_{<\C}$ the evident indexed full subcategory of $\D/(-)$ with value $(\D/I)_{<\C}$ at $I \in \D$.

Intriguingly, locally $\C$-small maps inherit closure under composition from $\C$.

 \begin{proposition}\label{prop:lscomp} 
    Locally small maps are closed under composition.
\end{proposition}

\begin{proof}
    Let $f : A \to B$ and $g : B \to C$ be locally small maps. Then, for any $D \in \C$ and $h : D \to C$, we can form successive pullbacks  
    \begin{equation*}
			\begin{tikzcd}
  (g^*h)^* A\ar[rr, "f^*g^*h"]\ar[dd, "(g^*h)^* f"'] \ar[ddrr, phantom, "\lrcorner", very near start]&& A\ar[dd, "f"]  \\
				\\
				h^*B \ar[rr, "g^*h" description]\ar[dd, "h^*g"'] \ar[ddrr, phantom, "\lrcorner", very near start]&& B\ar[dd, "g"]  \\
				\\
				D \ar[rr, "", "h"'] && C,
			\end{tikzcd}
		\end{equation*}
  in which $h^*B$ and $(g^*h)^* A$ are in $\C$ by the local smallness of $g$ and $f$, respectively. By the two pullbacks lemma, this $(g^*h)^* A\in \C$ is a pullback of $g \circ f$ along $h$. Thus, $g \circ f$ is locally small.
\end{proof}

We will be interested in the case where $\D$ is a topos and $\C$ is a logical subtopos. When $\F$ is a topos, we consider only  logical subtoposes of $\F$ that are full and closed under a choice of topos-theoretic structure in $\F$. In the rest of this section, we work in a topos $\F$ equipped with such a logical subtopos $\E$.

\begin{prop}\label{lem:logsubmon}
    All monomorphisms 
    are locally
    small maps.
\end{prop}

\begin{proof}
When $E \monoto B$ is a monomorphism of $\F$ and we have $A \in \E$ and $f : A \to B$, we form the diagram 
\begin{equation*}
			\begin{tikzcd}
				\{ a \in A ~|~ (\chi_E \circ f)(a)\} \ar[rrrr, bend left, "!"] \ar[rr, dashed]\ar[dd, rightarrowtail, ""'] \ar[ddrr, phantom, "\lrcorner", very near start]&& E\ar[rr, "!"]\ar[dd, rightarrowtail, ""'] \ar[ddrr, phantom, "\lrcorner", very near start]&& 1\ar[dd, rightarrowtail, "\top"]  \\
				\\
				A \ar[rr, "", "f"'] && B \ar[rr, "", "\chi_E"'] && \Omega
			\end{tikzcd}	\end{equation*}
in $\F$. By fullness, $\chi_m \circ f$ is a morphism of $\E$. Its comprehension $\{ a \in A ~|~ (\chi_E \circ f)(a)\}$ is thus in $\E$. As this is also a pullback of $E \monoto B$ along $f$, by the two pullbacks lemma, $E \monoto B$ is a locally small map.
\end{proof}

 In the present, topos-theoretic context, asking for a generic locally small map turns out to be equivalent to asking for a category $\pmb{\mho}$ in $\F$ subject to the desired pseudonatural equivalence  
 \begin{equation*}\label{eqn:unrep}
     \Hom_{\Cat(\F)}(I,\pmb{\mho}) \simeq (\F/I)_{<\E},
 \end{equation*}
where $I \in \F$. When $p : \mho_\bullet \to \mho$ is a generic locally small map, the corresponding $\pmb{\mho}$ is the internal full subcategory associated to $p$ \citep[B\'enabou 1970s; see][]{street1980cosmoi,street2018categories}. Since the notion of generic locally small map is closer to standard formulations \citep[\cf][]{benabou1973problemes,hyland1989theory,streicher2005universes}, we  identify universes in $\F$ with generic locally small maps.

However, our logical subtopos \[\E \inc \F\] is too weak a setting for the associated notion of generic locally small map to reflect set- or type-theoretic intuitions about universes. In light of prior approaches \citep{benabou1973problemes,hyland1989theory,streicher2005universes}, we certainly want the $\F$-indexed full subcategory
\[(\F/(-))_{<\E} \inc \F/(-)\]
to be an indexed logical subtopos. Indeed, this is the last ingredient needed for the  locally small maps to constitute an impredicative universe in the sense of \citet{streicher2005universes} (Definition \ref{def:imuniv}), which is needed to interpret the calculus of constructions.

We will thus adopt in Section \ref{sec:denslogsub} an additional density condition on our logical inclusion \[ \E \xhookrightarrow{i} \F.\] This will prove to be the natural way of extending the induced $\E$-indexed logical inclusion
\[\E/(-) \inc i^* \F/(-)\]
 to the desired $\F$-indexed logical inclusion
\[(\F/(-))_{<\E} \inc \F/(-).\]

\subsection{Dense Functors and Dense Subcategories}\label{sec:embden'}

This section houses category-theoretic preliminaries for our notion of dense logical subtopos (\S \ref{sec:denslogsub}), which we will use to axiomatize universes in a topos. These include an apparently novel definition of dense functor, which will be suggestive for us.

\subsubsection{Dense Functors}

Our account of dense functors \citep{ulmer1968properties} and related notions  will draw on sheaf- and stack-theoretic terminology and intuitions. The basic idea is that any functor $F : \C \to \D$ imbues each object $D$ of $\D$ with a canonical `cover' consisting of the comma category $F \downarrow D$, and, indeed, $F$ may be viewed as a `coverage' on $\D$. In this conception, a dense functor $F : \C \to \D$ is viewed as a `subcanonical coverage' on $\D$.  When $F$ is an inclusion of a full subcategory, this analogy may be cashed out by inducing an actual Grothendieck coverage on $\D$ from $F$, such that $F$ is dense \ifoif{} the Grothendieck coverage is subcanononical (See \S \ref{sec:denssub}).

When $F : \C \to \D$ is a  functor, we refer to the pair $(\D, F)$ as a \textbf{site} and say that $F$ is a \textbf{coverage} of $\D$. 

We recall that, when $\C$ and $\D$ are small categories, a functor 
\[F: \C \to \D\] 
induces the adjunction
\[F^* \dashv F_*: \PshC \to \PshD.\]
Here, $F^*$ is the precomposition functor, while $F_*$ is less intuitive. Nevertheless, we can readily derive a formula for $F_*$ in terms of $F^*$. For any $P \in \PshC$ and $I \in D$, we must have 
\begin{align*}
F_*P(I) & \cong \Hom_{\widehat{\D}}(\y I, F_* P) \\
	& \cong \Hom_{\widehat{\C}}( F^*(\y I),  P).
\end{align*}
For convenience, we officially set 
\[F_*P(I) : = \Hom_{\widehat{\C}}( F^*(\y I),  P).\]
When $P \in \PshD$ and $I \in \D$, we call an element of $F_*F^*P (I) = \Hom_{\widehat{\C}}( F^*(\y I),  F^* P)$ a \textbf{matching family} for $P$ at $I$. Such a matching family coherently assigns morphisms of form $f : FJ \to I$ to elements of $P(FJ)$. The unit 
\[\eta_P : P   \to F_*F^* P \]
at $P$ of the adjunction $F^* \dashv F_*$ maps an element of $P$ to its associated \textbf{canonical matching family}. For $I\in D$ and $e \in P(I)$, this is given by the assignment 
\[(f : FJ \to I) \mapsto P(f)(e).\]
We say that the element $e \in P(I)$ is an \textbf{amalgamation} for a matching family $\alpha$ for $P$ at $I$ if $\alpha$ is the canonical matching family associated to $e$.

The foregoing extends from presheaves to indexed categories.

Our 
functor \[F : \C \to \D\] also induces the pseudoadjunction of $2$-functors
\[F^* \dashv F_*: \IndCat(\C) \to \IndCat(\D)\]
 \citep[see][\S 3.2]{caramello2021relative}.
Here, $F^*$ is the precomposition $2$-functor. As for $F_*$, for any $P \in \IndCat(\C)$ and $I \in D$, we must have 
\begin{align*}
F_*P(I) & \simeq \Hom_{\IndCat(\D)}(\y I, F_* P) \\
	& \simeq \Hom_{\IndCat(\C)}( F^*(\y I),  P).
\end{align*}
For convenience, we officially set 
\[F_*P(I) : = \Hom_{\IndCat(\C)}( F^*(\y I),  P).\]
When $P \in \IndCat(\D)$ and $I \in \D$, we call an element of $F_*F^*P (I) =\Hom_{\IndCat(\C)}( F^*(\y I),  F^* P)$ a \textbf{descent datum} for $P$ at $I$. Such a descent datum consists of an assignment of morphisms of form $f : FJ \to I$ to objects of $P(FJ)$, together with appropriate coherences. The unit 
\[\eta_P : P   \to F_*F^* P \]
at $P$ of the pseudoadjunction $F^* \dashv F_*$ maps an object of $P$ to its associated \textbf{canonical descent datum}. For $I\in D$ and $e \in P(I)$, this is given by the assignment 
\[(f : FJ \to I) \mapsto P(f)(e),\]
together with the evident coherences. We say that the object $e \in P(I)$ is an \textbf{pseudoamalgamation} for a descent datum $\alpha$ for $P$ at $I$ if $\alpha$ is isomorphic to the canonical descent datum associated to $e$.

\begin{defn}
When $\C$ and $\D$ are small categories and $F : \C \to \D$ is a  functor, we say that a presheaf $P : \D\op \to \Set$ is a \textbf{sheaf} on $(\D,F)$ if the unit $\eta_P : P \to F_*F^* P$ at $P$ of the adjunction $F^* \dashv F_* : \PshC \to \PshD$ is an isomorphism. We also say that an object $D$ of $\D$ is a \textbf{sheaf} on $(\D,F)$ if the presheaf $\y D$ is a sheaf on $(\D,F)$. 
We say that a small $\D$-indexed category $P$ is a \textbf{stack} on $(\D,F)$ if the unit $\eta_P : P \to F_*F^* P$ at $P$ of the pseudoadjunction $F^* \dashv F_* : \IndCat(\C) \to \IndCat(\D)$ is an equivalence. 
We denote by $\Sh(\D,F)$ the full subcategory of $\PshD$ spanned by the sheaves on $(\D,F)$.  We denote by $\St(\D,F)$ the $2$-category of small stacks on $(\D,F)$, $\D$-indexed functors, and $\D$-indexed natural transformations.

\end{defn}

We are now ready to give our formulation of the notion of dense functor.

\begin{defn}
When $\C$ and $\D$ are small categories and $F : \C \to \D$ is a  functor, we say that the site $(\D,F)$ is \textbf{subcanonical} if every object $D \in \D$ is a sheaf on $(\D,F)$. In this case, we also say that the functor $F : \C \to \D$ is \textbf{dense}.
\end{defn}

This definition of dense functor agrees with the usual one:

\begin{prop}\label{prop:denseequiv}
Let $\C$ and $\D$ be small categories and $F : \C \to \D$  a  functor. Then, $F$ is dense \ifoif{} every object $D$ of $\D$ is the colimit of the diagram in $\D$ given by the composite of
\[F\downarrow D \xrightarrow{p_1} \C \xrightarrow{F} \D .\]
\end{prop}

\begin{proof}
($\Rightarrow$). Let $D$ and $X$ be objects of $\D$. We have
\begin{align*}
    \Hom_\D(D,X) & = \y X (D) \\
                 & \cong F_* F^* \y X (D) \\
                 & =  \Hom_{\widehat{\C}}(F^*\y D,  F^* \y X)\\
                  & =  \Hom_{\widehat{\C}}(\Hom_\D(F(-),D),\Hom_\D(F(-),X)).
\end{align*}
A bit of reflection reveals $\Hom_{\widehat{\C}}(\Hom_\D(F(-),D),\Hom_\D(F(-),X))$ as isomorphic to the set of cocones over the diagram $(F \circ p_1) : F\downarrow D \to \D$ with vertex $X$, and we conclude that $D$ is the colimit of this diagram.

($\Leftarrow$). Let $X$ and $D$ be objects of $\D$. Since we can write $D$ as a colimit of $F \circ p_1$, we have $\Hom_\D(D,X) \cong \Hom_{\widehat{\C}}(\Hom_\D(F(-),D),\Hom_\D(F(-),X))$, as before. We thus have 
\begin{align*}
    \y X (D) & = \Hom_\D(D,X) \\
    & \cong \Hom_{\widehat{\C}}(\Hom_\D(F(-),D),\Hom_\D(F(-),X)) \\
    & =  \Hom_{\widehat{\C}}(F^*\y D,  F^* \y X) \\
                 & = F_* F^* \y X (D)  .
\end{align*}
Thus, $X$ is a sheaf and $F$ is dense.
\end{proof}
 
\subsubsection{Dense Subcategories}\label{sec:denssub}

In practice, we will be interested in dense subcategories, rather than the more general dense functors. This will simplify matters, and we will see that it allows us to apply results from sheaf and stack theory. 

When  $i : \C \inc \D$ is an inclusion of a full subcategory, we write $(\D,\C)$ for $(\D, i)$  and say that $\C$ is a \textbf{coverage} of $\D$.

When $\C$ is  a full subcategory of $\D$, we denote by $J_\C$  the Grothendieck coverage of $\D$ defined by setting \[S \in J_\C (D) \Leftrightarrow \abs{\C \downarrow D} \subseteq S ,\]
for each object $D$ of  $\D$ and sieve $S$ on $D$. We call $J_\C$ (resp. $(\D, J_\C)$) the Grothendieck coverage (resp. site) \textbf{induced} by $\C$ (resp. $(\D, \C)$).

  It is straightforward to observe the following:

  \begin{prop}\label{cor:dincpsh}
Let $\C$ be  a full subcategory of a small category $\D$. Then, $\Sh(\D,\C) = \Sh(\D, J_\C)$ and $\St(\D,\C) = \St(\D, J_\C)$.
\end{prop}

We can thus apply results from sheaf and stack theory to our setting.

Though we will not use it, it is worth noting the following fact:

\begin{prop}\label{prop:functfac}
Let $\D$ be a small category and $i : \C \inc \D$ an inclusion of a full subcategory. Then, 
the composite
\[\Sh(\D,\C) \inc \PshD \xrightarrow{i^*} \PshC\]
is an equivalence.

\end{prop}

\begin{proof}
See \citet[Example C2.2.4(d)]{johnstone2002sketches}.
\end{proof}

By analogous arguments, one can show that also the composite
\[\St(\D,\C) \inc \IndCat(\D) \xrightarrow{i^*} \IndCat(\C)\]
is a biequivalence in this case.

\subsection{Dense Logical Subtoposes}\label{sec:denslogsub}

We are now in position to introduce our notion of \textbf{dense logical subtopos}.

\begin{defn}\label{def:denssub}
We say that a logical subtopos $\E$ of a small topos $\F$ is  \textbf{dense} if
\begin{enumerate}
    \item\label{ax:1} every object of $\F$ is a sheaf on $(\F,\E)$ (\ie{} $\E$ is a dense subcategory of $\F$);
    \item\label{ax:2} the self-indexing $\F/(-)$ is a stack on $(\F,\E)$.
\end{enumerate}
\end{defn}
Very roughly speaking, the novel axiom (\ref{ax:2}) requires that morphisms of $\F$ are determined by their reindexings to objects of $\E$. (One important caveat is that the reconstruction of a morphism of $\F$ from these reindexings additionally uses the coherence data that makes up the rest of a descent datum.)  The set-theoretic intution (see Example \ref{ex:set}) is that an arbitrary class function $f : X \to Y$ can be reconstructed from its reindexings to sets.

Axioms \ref{ax:1} and \ref{ax:2} assert that $\E$ has some control over $\F$, and may be regarded from an algebraic set theory perspective as `\textbf{projection principles}'  (converse but with similar motivation to  reflection principles, which would assert a similarity of $\F$ to $\E$). Our forthcoming use case for the notion of dense logical subtopos involves `projecting' properties like Cartesian closure from the $\mcE_m$-indexed category $\mcE_m/(-)$ to the $\mcE_n$-indexed category $(\mcE_n/(-))_{<\E_m}$..

The following examples demonstrate that dense logical subtoposes are natural-enough gadgets. We assume for the purpose of the examples a Grothendieck universe $\mathrm{Set}_0$, and denote by $\cat{Set}_0$ the full subcategory of $\Set$ spanned by the elements of $\mathrm{Set}_0$.

\begin{example}\label{ex:set}
    The category $\cat{Set}_0$ is a dense logical subtopos of $\Set$. Indeed, the terminal category $\cat{1} \inc \cat{Set}_0$ is already a dense subcategory of $\Set$ and $\Set/(-)$ is already a stack on $(\Set,\cat{1})$. The latter assertion is down to the fact that any function can be reconstructed from its fibers, which can be identified with reindexings to the terminal object $\cat{1}$.
\end{example}

\begin{example}\label{ex:psh}
   More generally, for any $\mathrm{Set}_0$-small category $\C$, the category of $\mathrm{Set}_0$-small presheaves $(\cat{Set}_0)^{\C\op}$ is a dense logical subtopos of the presheaf category $\Set^{\C\op}$. Indeed, the image of $\C$ under the Yoneda embedding $\y(\C) \inc (\cat{Set}_0)^{\C\op}$ is already a dense subcategory of $\Set^{\C\op}$ and $\Set^{\C\op}/(-)$ is already a stack on $(\Set^{\C\op},\y(\C))$. The latter assertion is down to the fact that any morphism of presheaves can be reconstructed from its fibers, together with appropriate functoriality data, which can collectively be identified with its canonical descent datum, when regarded as an object of the indexed category $\Set^{\C\op}/(-)$ over the site $(\Set^{\C\op},\y(\C))$.
\end{example}

\begin{example}\label{ex:sh}
  Even more generally, for any $\mathrm{Set}_0$-small site $(\C,J)$, the category of $\mathrm{Set}_0$-small sheaves $\Sh_0(\C,J)$ is a dense logical subtopos of the sheaf category $\Sh(\C,J)$. Indeed, the full image of $\y(\C)$ under sheafification $a(\y(\C)) \inc \Sh_0(\C,J)$ is already a dense subcategory of $\Sh(\C,J)$ and $\Sh(\C,J)/(-)$ is already a stack on $(\Sh(\C,J),a(\y(\C)))$. The latter assertion follows from the fact that $\Set^{\C\op}/(-)$ is a stack on $(\Set^{\C\op},\y(\C))$. 
  
  For any descent datum $\alpha$ at a sheaf $S$ for the indexed category $\Sh(\C,J)/(-)$ over the site $(\Sh(\C,J),a(\y(\C)))$, we can define a descent datum $\widetilde{\alpha}$ at $S$ for the stack $\Set^{\C\op}/(-)$ over the site $(\Set^{\C\op},\y(\C))$. This is given by the assignment
  \[(e : \y I \to S) \mapsto (\eta_{\y I})^* (\alpha_{a \y I} (\widetilde{e})),\]
  in which we use the unit and transposition for the adjunction $a \dashv i : \Set^{\C\op} \inc \Sh(\C,J)$, together with evident coherences. The pseudoamalgamation for $\widetilde{\alpha}$ is a morphism of presheaves $p: R \to S$. Since $R$ turns out to be a sheaf, $p$ also serves as a pseudoamalgamation for $\alpha$. 
\end{example}

Thus, the assumption of a Grothendieck universe yields an example of a dense logical subtopos in any (suitable) Grothendieck topos. For examples of `non-Grothendieck' dense logical subtoposes, it is natural to look to realizability toposes, though the picture there is not clear at present.

\begin{example}\label{ex:real}
   Given a $\mathrm{Set}_0$-small partial combinatory algebra $A$, we will write $\RT(A)$ (resp. $\PAss(A)$) for the realizability topos (resp. category of partitioned assemblies) associated to $A$. We will write $\RT_0(A)$ (resp. $\PAss_0(A)$) for the full subcategory of $\RT(A)$ (resp. $\PAss(A)$) associated to $A$ on objects with underlying set in $\mathrm{Set}_0$. 
   
   Is $\RT_0(A)$ a dense logical subtopos of $\RT(A)$? It is well-known that $\PAss(A)$ is a dense subcategory of $\RT(A)$, and $\PAss_0(A)$ is analogously a dense subcategory of $\RT_0(A)$. If $\PAss_0(A)$ was also a dense subcategory of $\RT(A)$, it would easily follow that $\RT_0(A)$ is a dense subcategory of $\RT(A)$, and we'd be off to a good start. 
   
   However, \citet[][proof of Theorem 1.5]{vanoosten2006filtered} suggests that this is not provable. More precisely, it suggests that $\mathrm{Set}_0$ would have to be (at least) a weakly compact cardinal to ensure that it does not admit any $\mathrm{Set}_0$-Aronszajn trees. The question of what (if any) large cardinal assumptions can be put on $\mathrm{Set}_0$ to prove that $\PAss_0(A)$ is a dense subcategory of $\RT(A)$ is open.
\end{example}

We now derive some important facts about dense logical subtoposes.

\begin{lemma}\label{prop:smstack}
    When $\E$ is a dense logical subtopos of a small topos $\F$, the $\F$-indexed category of locally $\E$-small maps $(\F/(-))_{<\E}$ is a stack on $(\F,\E)$.
\end{lemma}
\begin{proof}
  Since we have a commutative diagram
  \begin{equation*}
			\begin{tikzcd}
				 i^*((\F/(-))_{<\E})\ar[ddr, hook, ""'] \ar[rr, hookleftarrow, "\sim"] && \E/(-) \ar[ddl, hook]  \\
				\\
				& i^*(\F/(-))
			\end{tikzcd}
		\end{equation*}
  in $\IndCat(\E)$, we can build the commutative diagram
   \begin{equation*}
			\begin{tikzcd}
				(\F/(-))_{<\E} \ar[rr, "\eta_{(\F/(-))_{<\E}}"] \ar[dd, hook, ""]&& i_*i^*((\F/(-))_{<\E})\ar[dd, hook, ""'] \ar[r, hookleftarrow, "\sim"]& i_* (\E/(-)) \ar[ddl, hook]  \\
				\\
				\F/(-) \ar[rr, "\eta_{\F/(-)}"', "\sim"] && i_*i^*(\F/(-))
			\end{tikzcd}
		\end{equation*}
  in $\IndCat(\F)$. 
  It thus suffices to show that the restriction of 
  \[i_*i^*(\F/(-)) \xrightarrow{(\eta_{\F/(-)})\inv} \F/(-)\] 
  to $i_* (\E/(-))$ factors through $(\F/(-))_{<\E}$. But, indeed, any pseudoamalgamation in $\F/(-)$ for a descent datum of morphisms of $\E$ is locally $\E$-small by definition.
\end{proof}

\begin{proposition}\label{prop:denslog}
    When $\E$ is a dense logical subtopos of a small topos $\F$, $(\F/(-))_{<\E}$ is an indexed  logical subtopos of $\F/(-)$.
\end{proposition}

\begin{proof}
    In the proof of Lemma \ref{prop:smstack}, we established  a commutative diagram
    \begin{equation*}
			\begin{tikzcd}
				(\F/(-))_{<\E} \ar[rr, "\sim"] \ar[dd, hook, ""]&& i_*i^*((\F/(-))_{<\E})\ar[dd, hook, ""'] \ar[r, hookleftarrow, "\sim"]& i_* (\E/(-)) \ar[ddl, hook]  \\
				\\
				\F/(-) \ar[rr, "", "\sim"] && i_*i^*(\F/(-))
			\end{tikzcd}
		\end{equation*}
  in $\IndCat(\F)$. It thus suffices to show that $i_* (\E/(-))$ is an indexed logical subtopos of $i_*i^*(\F/(-))$. 

  For any $I \in \F$, the topos structure of $i_*i^*(\F/(-))(I)$ is computed pointwise using the indexed topos structure of $i^*(\F/(-))$. For example, the product $\alpha \times \beta$ of descent data $\alpha$ and $\beta$ for $\F/(-)$ at $I$ is given on object $J \in \E$ by the assignment 
  \[(f : J \to I) \mapsto \alpha_J(f) \times \beta_J(f),\]
in which $\times$ denotes the product in $\F/J$.
  
  Thus, since, as an indexed logical subtopos of $i^*(\F/(-))$, $\E/(-)$ is closed under the computation of topos structure in $i^*(\F/(-))$, $i_*(\E/(-))(I)$ is closed under the computation of topos structure in $i_*i^*(\F/(-))(I)$.
\end{proof}

\begin{corollary}\label{cor:streich}
 When $\E$ is a dense logical subtopos of a small topos $\F$, such that $\F$ admits a generic locally $\E$-small map, $(\F/(-))_{<\E}$ describes an impredicative universe in $\F$ in the sense of \citet{streicher2005universes}.
\end{corollary}

\begin{proof}
   The composition of two locally $\E$-small maps is a locally $\E$-small map, by Proposition \ref{prop:lscomp}. We have a generic locally $\E$-small map by assumption. Finally, $(\F/(-))_{<\E}$ is an indexed logical subtopos of $\F/(-)$, by  Proposition \ref{prop:denslog}. 
\end{proof}

\subsection{Stratified Toposes}\label{sec:strattop}

In this section, we introduce \textbf{stratified toposes}, though most of the legwork is behind us.  

Given a functor $\E: \N \to \cat{Topos}$, and  $m\leq n$, we may abusively write  

$ \E(m) \leq \E(n)$, for the evident functorial action of $\E$. If, moreover, $m < n$, we may abusively write 
$ \E(m) < \E(n)$, for $ \E(m) \leq \E(n)$.

\begin{defn}\label{def:nnnctm}
A (\textbf{small}) \textbf{stratified topos} consists of a functor \[\E: \N \to \cat{Topos}\]
such that 
\begin{itemize}
    \item $\E(m) \leq \E(n)$ is a dense logical subtopos inclusion, for all  $m \leq n$;
    \item  $\E(m) < \E(n)$ admits a generic locally 
    small map, for all  $m < n$.
\end{itemize}
\end{defn}

Henceforth, we consider only small stratified toposes without further comment.

\begin{rem}
  By Corollary \ref{cor:streich}, we see that, for any stratified topos $\E$ and  $m < n$, $(\E(n)/(-))_{<\E(m)}$ describes an impredicative universe in $\E(n)$ in the sense of \citet{streicher2005universes}.
\end{rem}

\begin{rem}
As discussed in Section \ref{sec:unibasic}, a conceptual selling point for the present approach is that we start from the notion of logical inclusion of toposes 
\[\E \inc \F,\]
rather than the less fundamental notion of indexed logical inclusion of form \[\mathbb{E} \inc \F/(-).\]

 Our approach now results in a  particularly nice characterization of a universe as an `untruncated subobject classifier.'
 
    In an ordinary topos $\E$, the subobject classifier $\pmb{\Omega}$, viewed as an internal poset, represents subobjects in the sense that \[\pmb{\sub} \cong [\pmb{\Omega}]\] 
    in $\Pos^{\E\op}$.
    
When $\E$ is a stratified topos, consider, for $m < n$, the internal full subcategory $\pmb{\mho}_m$ of $\E_n$ associated to a generic locally $\E_m$-small map $p_m : (\mho_\bullet)_m \to \mho_m$. Provided that the inclusion $i : \E_m \inc \E_n$ is not an equivalence (in which case $\E$ would be inconsistent), $\pmb{\mho}_m$ cannot represent arbitrary maps of $\E_m$ in the sense that 
\[\E_m/(-) \simeq [\pmb{\mho}_m];\]
 the internal category $\pmb{\mho}_m$ is too large to be in $\E_m$, so this equation doesn't type check.

We do have the next best thing, namely  \[i_*(\E_m/(-)) \simeq [\pmb{\mho}_m]\]
in $\IndCat(\E_n)$. For we have 
\begin{align*}
    i_*(\E_m/(-)) & \simeq i_*(i^*((\E_n/(-))_{<\E_m})) \\
    & \simeq (\E_n/(-))_{<\E_m} \\
    & \simeq [\pmb{\mho}_m] .
\end{align*}

By comparison, given an impredicative universe $\mathrm{E}$ in a topos $\F$ in the sense of \citet{streicher2005universes}, we have the analogous equation \[i_*(\mathbb{E}(1)) \simeq [\pmb{\mho}]\]
in $\IndCat(\F)$ (where $\pmb{\mho}$ is the internal full subcategory of $\F$ associated to a generic  $\mathrm{E}$-small map $p_m : \mho_\bullet \to \mho$) only only in the special case where $\mathbb{E}$ is a stack on $(\F,\mathbb{E}(1))$ (\cf{} Section \ref{sec:streicher}). 
\end{rem}

\begin{rem}
  Unfortunately, an analysis of the several differences between our stratified toposes and \cite{moerdijk2002type}'s stratified pseudotoposes would take us too far afield. However, we note that stratified pseudotoposes involve a nested sequence of pretoposes, rather than toposes, due to Moerdijk and Palmgren's predicative motivation. Our notion of stratified topos is geared towards modeling the calculus of constructions, with its impredicative type of propositions.   
\end{rem}

The following examples demonstrate that stratified toposes are natural-enough gadgets. We assume for the purpose of the examples a hierarchy of Grothendieck universes 
\[\mathrm{Set}_0, \mathrm{Set}_1, … \mathrm{Set}_n, …\]
indexed over $\N$, and denote by $\Set_n$ the full subcategory of $\Set$ spanned by the elements of $\mathrm{Set}_n$.

\begin{example}
    The sequence 
   \[\cat{Set}_0 \inc \cat{Set}_1 \inc … \cat{Set}_n \inc …\] 
   yields a stratified topos within $\Set$. For all  $m \leq n$, $\Set_m$ is a dense logical subtopos of $\Set_n$, following Example \ref{ex:set}. When, moreover, $m < n$, the projection function $p_m : (\mathrm{Set}_\bullet)_m \to \mathrm{Set}_m$ from pointed $\mathrm{Set}_m$-small sets is a generic locally  $\cat{Set}_m$-small map in $\cat{Set}_n$.
\end{example}

\begin{example}\label{ex:psh2}
 More generally, for any $\mathrm{Set}_0$-small category $\C$, the sequence 
   \[(\cat{Set}_0)^{\C\op} \inc (\cat{Set}_1)^{\C\op} \inc … (\cat{Set}_n)^{\C\op} \inc …\] 
   yields a stratified topos within $\PshC$. For all $m \leq n$, $(\cat{Set}_m)^{\C\op}$ is a dense logical subtopos of $(\cat{Set}_n)^{\C\op}$, following Example \ref{ex:psh}. When, moreover, $m < n$, a generic locally  $(\cat{Set}_m)^{\C\op}$-small map 
\[(\mho_m)_\bullet \to \mho_m\] 
in $(\cat{Set}_n)^{\C\op}$ is, following \citet{hofmann1999lifting}, given on $I \in \C$ by the projection function
\[\abs{1/((\cat{Set}_m)^{(\C/I)\op})} \to \abs{(\cat{Set}_m)^{(\C/I)\op}},\]
with the evident naturality in $I$.
\end{example}

\begin{example}
  Even more generally, for any $\mathrm{Set}_0$-small site $(\C,J)$, the sequence 
   \[\Sh_0(\C,J) \inc \Sh_1(\C,J) \inc … \Sh_n(\C,J) \inc …\] 
  yields a stratified topos within $\Sh(\C,J)$. For all  $m \leq n$, $ \Sh_m(\C,J)$ is a dense logical subtopos of $ \Sh_n(\C,J)$, following Example \ref{ex:sh}. When, moreover, $m < n$, a generic locally  $ \Sh_m(\C,J)$-small map in $\Sh_n(\C,J)$ may be obtained by, following \citet{streicher2005universes}, sheafifying the generic locally  $(\cat{Set}_m)^{\C\op}$-small map in $(\cat{Set}_n)^{\C\op}$ of Example \ref{ex:psh2}.
\end{example}

 Thus, the assumption of a tower of Grothendieck universes yields an example of a \textbf{stratified Grothendieck topos} in any (suitable) Grothendieck topos. By way of comparison, \citet{streicher2005universes} constructs his impredicative universes in any (suitable) Grothendieck or realizability topos. Due to the open problem discussed in Example \ref{ex:real}, we do not resolve the issue of whether a \textbf{stratified realizability topos} can be carved out of any (suitable) realizability topos. 

\section{The Stratified Topos of Coalgebras}\label{sec:univ'}

In this section, we get started on the study of stratified topos theory by constructing what we term the \textbf{stratified topos of coalgebras} for a \textbf{stratified Cartesian comonad} on a stratified topos.  
The proof is an interesting refinement of the classical construction of the topos of coalgebras for a Cartesian comonad, due to Lawvere and Tierney (1969 or 1970); universes in the stratified topos of coalgebras are constructed in much the same way as 
the subobject classifier in the topos of coalgebras.  
 
 \citet{moerdijk2002type} also raised, but left open, the issue of a coalgebra construction in the setting of stratified pseudotoposes;  our construction resolves the analogous issue in our comparable setting. A coalgebra construction for toposes equipped with an impredicative universe in the sense of \citet{streicher2005universes} was given in the dissertation of the author \citep[][Corollary 5.2.1.3]{zwanziger2023natural}, and provides a basis for the present construction. \citet{warren2007coalgebras} gave a comparable coalgebra construction for what he calls basic categories of classes. However, his construction of a generic display map in coalgebras is less natural; it requires much stronger assumptions than the approach of \citet{zwanziger2023natural}. It happens that these assumptions are not satisfied in Streicher's setting or the present one. Neither my construction nor Warren's applies to stratified pseudotoposes.

 Section \ref{sec:embden} extends our analysis of dense subcategories from Section \ref{sec:embden'}. Section \ref{sec:NMoUcomon} recalls and further develops  various aspects of the theory of comonads. Section \ref{sec:densecoalg} presents a coalgebra construction for dense logical subtoposes. Finally, Section \ref{sec:stratcoalg} constructs the stratified topos of coalgebras.

\subsection{Morphisms of Sites}\label{sec:embden}

This section houses category-theoretic preliminaries for our notion of stratified Cartesian functor. We extend the sheaf theory-inspired treatment of density from Section \ref{sec:embden'} with a corresponding notion of \textbf{morphism of sites}, which can perhaps be thought of as yielding a notion of morphism of dense subcategories.

\begin{defn}
Let $\D$ and $ \mathbf{F}$ be Cartesian categories, 
$ \C \inc \D$ and $\mathbf{E} \inc \mathbf{F}$ be Cartesian full subcategories, and \[G : \D \to \mathbf{F}\] a Cartesian functor. We say that $G$ is a \textbf{morphism of sites}, and write \[G : (\D,\C) \to (\mathbf{F}, \mathbf{E})\] if, for any object $D \in \D$ and morphism $f :  E \to GD \in \mathbf{F}$ such that $E \in \mathbf{E}$, there exists a morphism $g :  C \to D \in \D$ such that $C \in \C$ and $f$ factors through $Gg : GC \to GD$.

\end{defn}

It is straightforward to observe the following:

\begin{lem}\label{lem:togcov}
Let $\D$ and $ \mathbf{F}$ be Cartesian categories, 
$ \C \inc \D$ and $\mathbf{E} \inc \mathbf{F}$ be Cartesian full subcategories, and $G : \D \to \mathbf{F}$ a Cartesian functor. Then $G$ constitutes a morphism of sites  \[ (\D,\C) \to (\mathbf{F},\mathbf{E})\] \ifoif{}  it constitutes a morphism  \[ (\D,J_\C) \to (\mathbf{F},J_\mathbf{E})\] of induced Grothendieck sites.
\end{lem}

Lemma \ref{lem:togcov} allows us to use results from ordinary sheaf and stack theory:

\begin{lem}\label{lem:sheaf2sheaf}
Let $\D$ and $ \mathbf{F}$ be small Cartesian categories,  $ \C \inc \D$ and $ \mathbf{E} \inc \mathbf{F}$ be Cartesian full subcategories, and \[G : (\D, \C) \to (\mathbf{F},\mathbf{E})\]  a morphism of sites. Then, \[G^* : \Set^{\mathbf{F}\op} \to \PshD\] maps sheaves on $(\mathbf{F},\mathbf{E})$ to sheaves on $(\D, \C)$. Moreover, \[G^* : \IndCat(\mathbf{F})  \to \IndCat(\mathbf{D})\] maps stacks on $(\mathbf{F},\mathbf{E})$ to stacks on $(\D, \C)$.
\end{lem}

\begin{proof}
In light of Corollary \ref{cor:dincpsh} and Lemma \ref{lem:togcov}, we may apply statements from the literature about morphisms of Grothendieck sites. The relevant statements appear as \citet[][Lemma C2.3.3]
{johnstone2002sketches}  and \citet[][Proposition 3.4.2]{caramello2021relative}.

\end{proof}

\begin{lem}\label{lem:presls}
Let $\D$ and $ \mathbf{F}$ be Cartesian categories, 
$ \C \inc \D$ and $\mathbf{E} \inc \mathbf{F}$ be Cartesian full subcategories, and 
\begin{equation*}
			\begin{tikzcd}
				\D \ar[rr, "G",""']&& \mathbf{F}   \\
				\\
				\C \ar[uu, hook, ""]\ar[rr, "", "F"'] && \mathbf{E} \ar[uu, hook, ""']
			\end{tikzcd}
		\end{equation*}
a commutative square of Cartesian functors such that \[G : \D \to \mathbf{F}\] constitutes a morphism of sites \[G : (\D,\C) \to (\mathbf{F},\mathbf{E}).\] Then, $G$ takes locally $\C$-small maps to locally $\mathbf{E}$-small maps.
\end{lem}

\begin{proof}
    Let $p : E \to B$ be a locally $\C$-small map. Then, any morphism \[A \xrightarrow{f} G B\] with $A \in \C$ factors as  \[A \xrightarrow{h} GC \xrightarrow{Gg} GB,\] in which $C \in \C$. Consequently, we can form a pullback 
\begin{equation*}
			\begin{tikzcd}
				D \ar[dd, ""']\ar[ddrr, phantom, "\lrcorner", very near start]\ar[rr, ""]&& GE \ar[dd, "G p"]\\
				\\
				A  \ar[rr, "f"'] &&GB
			\end{tikzcd}
		\end{equation*}
  as the sequential pullback
  \begin{equation*}
			\begin{tikzcd}
				h^* (G(g^*E)) \ar[dd, ""']\ar[ddrr, phantom, "\lrcorner", very near start]\ar[rr, ""]&& G(g^*E) \ar[ddrr, phantom, "\lrcorner",very near start]\ar[rr, ""]\ar[dd, ""']&& GE \ar[dd, "G p"]\\
				\\
				A  \ar[rr, "h"']&& GC\ar[rr, "Gg"'] &&GB.
			\end{tikzcd}
		\end{equation*}
  Since $p$ is locally $\C$-small, we have $g^* E \in \C$, $G(g^*E) = F(g^*E) \in \mathbf{E}$, and, finally, $h^* (G(g^*E)) \in \mathbf{E}$.

\end{proof}

\begin{prop}\label{prop:mordensecol}
Let $\C,\D, \mathbf{E},$ and $ \mathbf{F}$ be small Cartesian categories, $ i: \C \inc \D$ and $ j: \mathbf{E} \inc \mathbf{F}$ Cartesian full subcategory inclusions, such that $j$ is dense, and 
\begin{equation*}
			\begin{tikzcd}
				\D \ar[rr, "G",""']&& \mathbf{F}   \\
				\\
				\C \ar[uu, hook, "i"]\ar[rr, "", "F"'] && \mathbf{E} \ar[uu, hook, "j"']
			\end{tikzcd}
		\end{equation*}
a commutative square of Cartesian functors such that \[G : \D \to \mathbf{F}\] constitutes a morphism of sites \[G : (\D,\C) \to (\mathbf{F},\mathbf{E}).\]

Then, for any $D \in \D$, $GD \in \mathbf{F}$ may be written as \[\varinjlim_{(C \in \C, f : C \to D)} G C.\]
\end{prop}

\begin{proof}
Let $X \in \mathbf{F}$. We have
\begin{align*}
    \Hom_\mathbf{F}(GD,X) & \cong \Hom_{\widehat{\mathbf{E}}}(\Hom_\mathbf{F}(j(-),GD),\Hom_\mathbf{F}(j(-),X)) & (\mbox{Proposition \ref{prop:denseequiv}})\\
     & \cong \Hom_{\widehat{\mathbf{C}}}(\Hom_\mathbf{D}(i(-),D), \Hom_\mathbf{F}(Gi(-), X) ) ,
\end{align*}
in which the last isomorphism follows from a diagram chase, notably by applying that $G$ is a morphism of sites. A bit of reflection on the resulting equation \[\Hom_\mathbf{F}(GD,X) \cong \Hom_{\widehat{\mathbf{C}}}(\Hom_\mathbf{D}(i(-),D), \Hom_\mathbf{F}(Gi(-), X) )\] reveals $GD$ as the desired colimit.
\end{proof}

\subsection{Comonads}\label{sec:NMoUcomon}

We will intensively use both the classical \citep{godement1958topologie} category-theoretic notion of a comonad $\f : \C \to \C$ on a (small) category $\C$ and the generalization of this from $\Cat$ to an arbitrary $2$-category \citep[see][]{street1972formal}.

Section \ref{sec:formth} reviews the so-called formal theory of comonads in a $2$-category. Section \ref{sec:appl} collects various applications of the theory to categories, indexed categories, and stacks. Section \ref{sec:dispcomon} is the most novel and deals with comonads in the context of display maps and locally small maps in particular.

\subsubsection{The Formal Theory of Comonads}\label{sec:formth}

In this section, we collect relevant information pertaining to the application of the formal theory of comonads as described by \citet{street1972formal}.

In this section, we work in the context of a $2$-category $\C$, considering only comonads, adjunctions, \etc, in this $2$-category.

		We frequently write $\f: C \to C$ for a comonad on an object $C$ with underlying morphism $\f : C \to C$, counit $\ve : \f \Rightarrow \id_C$ and comultiplication $\delta : \f \Rightarrow \f \f $. We frequently write $L \dashv R : C \to D $ for an adjunction with right and left adjoints $R : C \leftrightarrows D :L$, unit $\eta: \id_D \Rightarrow RL$, and counit $\ve : LR \Rightarrow \id_C$.

  We recall the following basic fact.
	
	\begin{prop}\label{indcom}
	Any adjunction $L \dashv R : C \to D$ induces a comonad $\f: C \to C$ by setting
		\begin{itemize}
			\item $\f :\equiv LR$;
			\item $\ve : \equiv \ve$;
			\item $\delta : \equiv L \eta R$.
		\end{itemize}
			\end{prop}

		Any adjunction $L \dashv R : C \to D$ that induces (by Proposition \ref{indcom}) the comonad $\f : C \to C$ is said to be a \textbf{decomposition} for $\f$.

When $\f_1 : C_1 \to C_1$ and $\f_2 : C_2 \to C_2$ are comonads,  a \textbf{morphism of comonads}, $f: \f_1 \to \f_2$,   consists of a morphism $f: C_1 \to C_2$ between the underlying objects, together with a $2$-morphism $\tau: f\f_1 \Rightarrow \f_2 f$ satisfying the commutative diagrams    
		\begin{center}
			\begin{tikzcd}
				f\f_1 \ar[rrrr, "\tau"] \ar[ddrr, "f \ve_1"']&&&& \f_2 f \ar[ddll, " \ve_2 f"]\\
				\\
				&&f
			\end{tikzcd}
		\end{center}
		and 
		\begin{center}
			\begin{tikzcd}
				f \f_1 \ar[rrrr, "\tau"]\ar[dd, "f \delta_1"'] &&&& \f_2  f \ar[dd, "\delta_2 f"] \\
				\\
				f \f_1\f_1  \ar[rr, "\tau \f_1"'] &&\f_2 f\f_1 \ar[rr, "\f_2 \tau"'] && \f_2 \f_2 f.
			\end{tikzcd}
		\end{center}

  We say that such a morphism of comonads is \textbf{strong} if $\tau$ is an isomorphism and \textbf{strict} if $\tau$ is an identity (in which case it is superfluous data).

		When $\f_1 : C_1 \to C_1$ and $\f_2 : C_2 \to C_2$ are comonads and $f_1 : \f_1 \to \f_2$ and $f_2 : \f_1 \to \f_2$ are morphisms of comonads, a \textbf{$2$-morphism of comonads}, $\alpha : f_1 \Rightarrow f_2 (: \f_1 \to \f_2)$, consists of a $2$-morphism, denoted $\alpha : f_1 \Rightarrow f_2 (:  C_1 \to C_2)$, between the underlying morphisms, satisfying the commutative diagram
		
		\begin{center}
			\begin{tikzcd}
				f_1 \f_1 \ar[dd, "\tau_1"'] \ar[rr, "\alpha \f_1"]&& f_2\f_2 \ar[dd, "\tau_2"]\\
				\\
				\f_2 f_1  \ar[rr, "\f_2 \alpha"'] &&\f_2 f_2.
			\end{tikzcd}
		\end{center}

	We write $\Comon(\C)$ for the $2$-category of comonads and morphisms and $2$-morphisms of comonads in $\C$. 

When $\f : C \to C$ is a comonad, an \textbf{object of coalgebras} for $\f$ consists of an object $C^\f$, such that 
\begin{equation}\label{eqn:pre1}
    \Hom_{\C}(D,C^\f) \cong \Hom_{\Comon(\C)}(\Id_D,\f)~,
\end{equation}
$2$-naturally in $D$.\footnote{We follow \citet{street1972formal} in using the term `object of coalgebras.' However, the less suggestive term `(co)Eilenberg-Moore object' is more common.} We say that $\C$ \textbf{admits the construction of coalgebras} if an object of coalgebras exists for each comonad in it. 

The assignment \[C \mapsto (\id_C : C \to C)\] extends to a $2$-functor 
\[\Inc : \C \to \Comon(\C)\] in the evident fashion. Our $\C$ admits the construction of coalgebras \ifoif{} $\Inc$ admits a right $2$-adjoint \[\Coalg : \Comon(\C) \to \C,\] due to general considerations about $2$-adjoints. On objects, $\Coalg$ is the assignment  \[(\Box : C \to C) \mapsto C^\Box.\]

When $C^\f$ is an object of coalgebras for a comonad $\Box : C \to C$, the \textbf{cofree morphism} \[F : C \to C^\f\] is defined as the transpose across the isomorphism \ref{eqn:pre1} of the morphism of comonads \[(\f, \delta) : (C, \id) \to (C, \f),\] while the \textbf{forgetful morphism} \[U : C^\f \to C\] is defined as the underlying morphism of the transpose across the isomorphism \ref{eqn:pre1} of the morphism \[\Id_{C^\f} : C^\f \to C^\f.\]
As developed in \citet{street1972formal}, we then have 
\[U \dashv F : C \to C^\f,\]
and this adjunction is a decomposition for \[\Box : C \to C,\]
which we term the \textbf{forgetful-cofree decomposition}.

\subsubsection{Assorted Applications}\label{sec:appl}

We collect here various general statements about comonads that we will use, which are less abstract than the theory presented in Section \ref{sec:formth}.

When $\C$ and $\D$ are categories and $\f_\C: \C \to \C$ and $\f_\D : \D \to \D$ are comonads, we say that a strict morphism of comonads \[ \f_\C \xrightarrow{i} \f_\D\]  is a \textbf{full inclusion} of comonads, and write \[\f_\C \xhookrightarrow{i} \f_\D,\] if the underlying functor \[ \C \xrightarrow{i} \D\] is a full subcategory inclusion.  When $\C$, $\D$, $\f_\C$, and $\f_\D$ are Cartesian, we say that a full inclusion of comonads \[ \f_\C \xhookrightarrow{i} \f_\D\]  is a \textbf{Cartesian full inclusion} of Cartesian comonads if the underlying functor \[ \C \xhookrightarrow{i} \D\] is Cartesian.

The following should be fairly evident:

\begin{lem}\label{lem:comoninclinc1}\label{lem:comoninclinc}
Let $\C$ and $\D$ be categories (resp. Cartesian categories), $\f_\C: \C \to \C$ and $\f_\D : \D \to \D$ be comonads (resp. Cartesian comonads), and \[ \f_\C \xhookrightarrow{i} \f_\D\] be a full  inclusion of comonads (resp. Cartesian full  inclusion of Cartesian comonads). Then, the functor \[\Coalg(\f_\C) \xrightarrow{\Coalg(i)} \Coalg(\f_\D)\] is a full subcategory inclusion (resp. Cartesian full subcategory inclusion). 
\end{lem}

When $\E$ and $\F$ are toposes and $\f_\E: \E \to \E$ and $\f_\F : \F \to \F$ are Cartesian comonads,  we say that a Cartesian full inclusion of Cartesian comonads  \[\f_\E \xhookrightarrow{i} \f_\F\]  
is a \textbf{logical inclusion}  of Cartesian comonads on toposes
if the underlying functor \[\E \xhookrightarrow{i}  \F\] is logical.

\begin{lem}\label{lem:comonloglog}
Let $\E$ and $\F$ be toposes, $\f_\E: \E \to \E$ and $\f_\F : \F \to \F$  be Cartesian comonads, and  \[\f_\E \xhookrightarrow{i} \f_\F\] a  logical inclusion of Cartesian comonads on toposes. Then, the Cartesian full subcategory inclusion \[ \Coalg(\f_\E) \xhookrightarrow{\Coalg(i)} \Coalg(\f_\F)\] is a  logical inclusion of toposes.
\end{lem}

\begin{proof}
The underlying functor $i : \E \inc \F$ preserves and commutes with everything in sight. As the reader may verify, this includes the construction of finite limits, exponentials, and the subobject classifier in coalgebras. The functor $\Coalg(i) : \Coalg(\f_\E) \to \Coalg(\f_\F)$ is thus logical. 
\end{proof}

When $\C$ is a small category and $\Box : P \to P$ is a comonad in $\IndCat(\C)$, the \textbf{indexed category of coalgebras} is the evident $\C$-indexed category $P^\Box$ given on objects $I\in \C$ by $P^\f(I) = P(I)^{(\f_I)}$. 

\begin{prop}\label{mainthpfib}
Given a small category $\C$, the $2$-category $\IndCat(\C)$ admits the construction of coalgebras.
\end{prop}

\begin{proof}
	Given a comonad $\f :  P \to P$ in $\IndCat(\C)$, the object of coalgebras is given by the indexed category of coalgebras $P^\f$.
\end{proof}

\begin{prop}\label{mainthpst}
Given a small site $(\C,J)$, the $2$-category $\St(\C,J)$ admits the construction of coalgebras.
\end{prop}

\begin{proof}
	$\St(\C,J)$ is closed under the construction of bilimits in $\IndCat(\C)$. Since objects of coalgebras are flexible limits \citep{bird1989flexible}, they are bilimits. Thus, given a comonad $\f :  P \to P$ in $\St(\C,J)$, the object of coalgebras is again given by the indexed category of coalgebras $P^\f$.
\end{proof}

\subsubsection{Display Comonads}\label{sec:dispcomon}

We collect in this section some critical development for our construction of the stratified topos of coalgebras, which may be unified via the use of display maps (see Section \ref{sec:streicher}).

When $\C$ is a category equipped with a class of display maps, we say that a comonad $\Box : \C \to \C$ is  \textbf{display} when the underlying functor
\begin{itemize}
    \item takes display maps to display maps, and 
    \item takes pullbacks of display maps to pullbacks of display maps.
 \end{itemize}

When $\C$ is a category and $\mathrm{D}$ is a class of display maps in $\C$, for which a comonad $\Box : \C \to \C$ is display, we can define a class of display maps $\mathrm{D}^\Box$ in $\C^\Box$ by  setting \[p \in \mathrm{D}^\Box \Leftrightarrow Up \in \mathrm{D}.\]

Moreover, we can define a $\C^\Box$-indexed comonad 
\[\mathbb{B} : U^*\mathbb{D} \to U^*\mathbb{D},\]
which we term the \textbf{induced indexed comonad}. 
 On $(B, b) \in \C^\f$ and a display map of form $p : A \to B$, this is given by 
 \[\mathbb{B}_{(B, b)}(p) = b^*(\Box (p)) .\]

 Similarly, the cofree functor $F : \C \to \C^\Box$ induces a $\C^\Box$-indexed functor  \[\mathbb{F} : U^* \mathbb{D} \to \mathbb{D^\Box},\]
 in which $\mathbb{D^\Box}$ is the indexed category of display maps corresponding to $D^\Box$.
 On $B \in \C^\Box$ and a display map of form $p : A \to UB$, this is given by 
 \[\mathbb{F}_{B}(p) = (\eta_B)^*(F (p)) .\] We term $\mathbb{F}$ the  \textbf{induced indexed right adjoint}, as it admits a $\C^\Box$-indexed left adjoint \[\mathbb{U} :  \mathbb{D}^\Box \to U^* \mathbb{D},\]
 given on $B \in \C^\Box$ and a display map $p : A \to B$ by 
 \[\mathbb{U}_{B}(p) = U(p).\] Of course, this \[\mathbb{U} \dashv \mathbb{F} : U^* \mathbb{D} \to \mathbb{D^\Box}\] we term the \textbf{induced indexed adjunction}. As the reader may verify, the induced indexed adjunction is a decomposition for the induced indexed comonad. This decomposition is comonadic in a suitable sense, though we will use only the following relation between the class of display maps $\mathrm{D^\Box}$ in $\C^\Box$ and the induced $\C^\Box$-indexed comonad $\mathbb{B} : U^*\mathbb{D} \to U^*\mathbb{D}$:

 \begin{proposition}\label{prop:comonadic}
     $\mathbb{D^\Box} \cong (U^*\mathbb{D})^{\mathbb{B}} \in \IndCat(\C^\Box)$.
 \end{proposition}

 \begin{proof}
      Given a $\Box$-coalgebra $(B, b) \in \C^\f$ and display map $p : E \to B$ in $\C$, we observe that $\Box$-coalgebra structures on $E$ making $p$ into a morphism of $\Box$-coalgebras (equivalently objects of $\mathbb{D^\f}(B, b)$ laying over $p$) are in correspondence with $\mathbb{B}_{(B,b)}$-coalgebra structures on $p$ (\ie{} objects of $(U^*\mathbb{D})^\mathbb{B}(B,b)$ laying over $p$), mediated by the diagram
 \begin{center}
	    	\begin{tikzcd}
	    	E \ar[dddr, "p"', bend right] \ar[rrrd, "", bend left] \ar[dr, dashed]\\
			&	\mathbb{B}_{(B,b)} E \ar[dd,"\mathbb{B}_{(B,b)} p" description]\ar[ddrr, phantom, "\lrcorner", very near start]\ar[rr, "" ]&&  \Box E\ar[dd, "\Box p"]\\
				\\
			&	B \ar[rr, "b"'] &&\Box B
			\end{tikzcd}
	\end{center}
in $\C$.  See also \citet[][Proposition 4.2.3]{zwanziger2023natural}.
 \end{proof}

 \begin{example}\label{ex:cartc}
     When $\C$ is a Cartesian category, we may regard its class of arrows $C_1$ as a class of display maps in $\C$. A Cartesian comonad \[\Box : \C \to \C\] is then display and induces an indexed comonad of form \[\mathbb{B} : \C/(U(-)) \to \C/(U(-)),\]
     and Proposition \ref{prop:comonadic} yields the identity \[\C^\Box/(-) \cong (\C/(U(-)))^\mathbb{B}.\]
 \end{example}

\begin{example}\label{ex:locsmall}
    When  $\D$ is a Cartesian category and $\C \inc \D$ is a Cartesian full subcategory, we may regard the class $(D_1)_{<\C}$ of locally $\C$-small maps in $\D$ as a class of display maps in $\D$. When $\f_\C: \C \to \C$ and $\f_\D : \D \to \D$ are Cartesian comonads, the inclusion \[\C \xhookrightarrow{i} \D\] 
 constitutes a Cartesian full inclusion of Cartesian comonads
 \[\f_\C \xhookrightarrow{i} \f_\D,\] and the underlying Cartesian functor 
\[ \D \xrightarrow{\f_\D} \D\]
constitutes a morphism of sites
\[ (\D,\C) \xrightarrow{\f_\D} (\D,\C),\] 
 $\f_\D$ takes locally $\C$-small maps to locally $\C$-small maps (Lemma \ref{lem:presls}). The comonad $\f_\D$ is  then display and induces an indexed comonad of form \[\mathbb{B} : (\D/U(-))_{<\C} \to (\D/U(-))_{<\C}.\]

 We note that, when $\C$ contains a subobject classifier 
 for $\D$, the Cartesian functor \[U_\D : \D^{\f_\D}  \to \D\]
constitutes a morphism of sites
\[U_\D : (\D^{\f_\D},\C^{\f_\C})  \to (\D,\C).\]

For, let $A \in \C$, $B \in \D^{\f_\D}$, 
and $f :  A \to U_\D B \in \D$. Then, because the Cartesian functor 
\[\f_\D : \D \to \D\]
constitutes a morphism of sites
\[\f_\D : (\D,\C)  \to (\D,\C),\]
there exists some 
$A' \in \C$, $g : A' \to U_\D B \in \D$, and  
$h : A \to \f_\D A' (=\f_\C A') \in \C$ such that the outside of the diagram
\begin{center}
\begin{tikzcd}
A \ar[rrrd,"h",bend left] \ar[dddr,"f"',bend right]\ar[dr, dashed, "(f{,}h)" description]\\
&\mathbb{B}_B A' \ar[dd, "\mathbb{B}_B g" description]\ar[ddrr, phantom, "\lrcorner", very near start]\ar[rr, rightarrowtail, "(\f_\D g)^*\ko_B" description]&& \f_\C A' \ar[dd, "\f_\D g"]\\
\\
&U_\D B \ar[rr, rightarrowtail, "\ko_B"'] &&\f_\D U_\D B
\end{tikzcd}
\end{center}
commutes in $\D$, in which $k_B$ denotes the coalgebra structure map of $B$. We may form the pullback shown, and obtain the morphism $(f,h) : A \to \mathbb{B}_B A'$ in $ \D$ from its universal property. 

This $\mathbb{B}_B A'$ is in $\C$: since \[ U_\D B \xrightarrow{\ko_B} \f_\D U_\D B\] is a coalgebra structure map, it is a split mono. In particular, it is a  
mono and thus locally $\C$-small (Lemma \ref{lem:logsubmon}). Its pullback along $\f_\D g$,  \[ \mathbb{B}_B A' \xrightarrow{(\f_\D g)^*\ko_B} \f_\D A' (= \f_\C A') ,\] then has $\mathbb{B}_B A'$ chosen to be in $\C$ because $\f_\C A'$ is in $\C$.

The morphism \[ \mathbb{B}_B A' \xrightarrow{\mathbb{B}_B g} U_\D B\] of $\D$ underlies the  morphism \[\mathbb{F}_B A' \xrightarrow{\mathbb{F}_B g} B\] of $\D^{\f_\D}$.  
This $\mathbb{F}_B A'$ is in $\C^{\f_\C}$ because its underlying object $ \mathbb{B}_B A'$ is in $\C$.

One consequence of $U_\D : \D^{\f_\D}  \to \D$ being a morphism of sites is that the morphisms of $\D^{\f_\D}$ with an underlying locally $\C$-small map are precisely the locally $\C^{\f_\C}$-small maps. When a morphism $p : E \to B$ of $\D^{\f_\D}$ has an underlying locally $\C$-small map, we have, for any $A \in \C^{\f_\C}$ and $f : A \to B$,  $U_\D(f^* E) = (U_\D f)^* (U_\D E) \in \C$. For the forgetful functor $U_\D : \D^{\f_\D} \to \D$ not only creates pullbacks, but also preserves them on the nose, provided that those in $\D^{\f_\D}$ are constructed using the choice of pullbacks in $\D$. Our $f^* E$ is thus in  $\C^{\f_\C}$ because its underlying object is in $\C$. So $p$ is a locally $\C^{\f_\C}$-small map.  Conversely, the morphism of sites  
\[U_\D : (\D^{\f_\D},\C^{\f_\C}) \to (\D,\C)\] takes any locally $\C^{\f_\C}$-small map to an underlying locally $\C$-small map (Lemma \ref{lem:presls} again).

Combining this last observation with Proposition \ref{prop:comonadic} yields the identity \[(\D^{\Box_\D}/(-))_{<\C^{\Box_\C}} \cong ((\D/U(-))_{<\C})^\mathbb{B}.\]

\end{example}

\subsection{The Dense Logical Subtopos of Coalgebras}\label{sec:densecoalg}

We now show that dense logical subtoposes admit a coalgebra construction, a key lemma for our construction of the stratified topos of coalgebras.

When $\E$ and $\F$ are toposes and $\f_\E: \E \to \E$ and $\f_\F : \F \to \F$ are Cartesian comonads, we say that a logical inclusion of Cartesian comonads on toposes \[\f_\E \inc \f_\F\] is  \textbf{dense} if its underlying logical inclusion of toposes \[\E \inc \F\] is dense and the underlying Cartesian functor 
\[ \F \xrightarrow{\f_\F} \F\]
constitutes a morphism of sites
\[ (\F,\E) \xrightarrow{\f_\F} (\F,\E).\]

\begin{lem}\label{lem:topdensedense}
Let $\E$ and $\F$ be toposes, $\f_\E: \E \to \E$ and $\f_\F : \F \to \F$ be Cartesian comonads, and  \[ \f_\E \xhookrightarrow{i} \f_\F\] be a  dense logical inclusion of Cartesian comonads on toposes. Then, 
the logical inclusion of toposes \[ \Coalg(\f_\E) \xhookrightarrow{\Coalg(i)} \Coalg(\f_\F)\] is dense.
\end{lem}

\begin{proof}

To show the density of $\Coalg(i) : \Coalg(\f_\E) \inc \Coalg(\f_\F)$ as a functor, it suffices (Proposition \ref{prop:denseequiv}) to show that each $\Delta \in \Coalg(\f_\F)$ is the colimit of the composite
\[\Coalg(i)\downarrow \Delta \xrightarrow{p_1} \Coalg(\f_\E) \xhookrightarrow{\Coalg(i)} \Coalg(\f_\F) . \]
Since comonadic functors create all colimits that exist in their codomain, it moreover suffices to show that $U_\F\Delta \in \F $  is the colimit of the composite
\[\Coalg(i)\downarrow \Delta \xrightarrow{p_1} \Coalg(\f_\E) \xhookrightarrow{\Coalg(i)} \Coalg(\f_\F) \xrightarrow{U_\F} \F .\]

In view of the  density of $i : \E \inc \F$ and 
the fact (Example \ref{ex:locsmall}) that that the functor \[U_\F : \Coalg(\f_\F)  \to \F\]
constitutes a morphism of sites
\[U_\F : (\Coalg(\f_\F),\Coalg(\f_\E))  \to (\F,\E),\]
this follows from Proposition \ref{prop:mordensecol}.

There remains to show that $\Coalg(\f_\F)/(-)$ is a stack on $(\Coalg(\f_\F), \Coalg(\f_\E))$. Again because  

the functor \[ \Coalg(\f_\F)  \xrightarrow{U_\F} \F\]
constitutes a morphism of sites
\[U_\F : (\Coalg(\f_\F),\Coalg(\f_\E))   \xrightarrow{U_\F} (\F,\E),\]
 $(U_\F)^*$ takes stacks on $(\F,\E)$ to stacks on $(\Coalg(\f_\F),\Coalg(\f_\E))$ (Lemma \ref{lem:sheaf2sheaf}). Thus,  $(U_\F)^* (\F/(-))$ 
is a stack on $(\Coalg(\f_\F),\Coalg(\f_\E))$. Since the inclusion $\St(\Coalg(\f_\F),\Coalg(\f_\E)) \inc \IndCat(\Coalg(\f_\F))$ is full on both morphisms and $2$-morphisms,  the induced indexed comonad $\mathbb{B} : (U_\F)^* (\F/(-)) \to (U_\F)^* (\F/(-))$ in $\IndCat(\Coalg(\f_\F))$  (derived from Example \ref{ex:cartc}) lies in $\St(\Coalg(\f_\F),\Coalg(\f_\E))$. By the proof of Proposition \ref{mainthpst},  $\Coalg(\mathbb{B})$ is then a stack on $(\Coalg(\f_\F),\Coalg(\f_\E))$. And  $\Coalg(\f_\F)/(-) \cong \Coalg(\mathbb{B})$, by Example \ref{ex:cartc}.
\end{proof}

\subsection{The Stratified Topos of Coalgebras}\label{sec:stratcoalg}

We are finally ready to construct the stratified topos of coalgebras.

For intuition, let us recall some of the classical construction of the topos of coalgebras (for details, see, \eg{}, \citealt{mac92}, \S V.8.). Let $\f : \mathcal{E} \to \mathcal{E}$ be a Cartesian comonad on a topos $\mathcal{E}$. We construct the category of coalgebras $\mathcal{E}^\f$ and the forgetful-cofree decomposition $U \dashv F: \mathcal{E} \to \mathcal{E}^\f $. We then verify that $\mathcal{E}^\f$ admits a subobject classifier by constructing one from a subobject classifier $\pmb{\Omega}_\mathcal{E} $ of $\mathcal{E}$. Curiously, the  comonad 
 \begin{equation}\label{eqn:51}
     \f : \mathcal{E} \to \mathcal{E}
 \end{equation}
 induces a comonad 
 \begin{equation}\label{eqn:52}
     \beta : F\pmb{\Omega}_{\mathcal{E}} \to F\pmb{\Omega}_{\mathcal{E}}
 \end{equation}
 of posets	\textit{in} $\mathcal{E}^\f$. 
	
	A subobject classifier $ \pmb{\Omega}_{\E^\f}$ for $\E^\f$ is then constructed as the internal poset of coalgebras (or, equivalently, the internal poset of fixed points) $(F \pmb{\Omega}_{\mathcal{E}})^\beta$. 
	
	We will undertake an untruncated version of the same construction for universes, working with internal categories rather than internal posets. We will let $\f : \mathcal{E} \to \mathcal{E}$ be what we call a stratified Cartesian comonad on a stratified topos $\mathcal{E}$. We will construct the corresponding object of coalgebras $\mathcal{E}^\f$ and the forgetful-cofree decomposition $U \dashv F : \mathcal{E} \to \mathcal{E}^\f $ in $\Topos^\N$. We will then verify that, for all $m<n$, $\mathcal{E}^\f(n)$ admits a universe for locally $\mathcal{E}^\f(m)$-small maps by constructing one from the universe $\pmb{\mho}_{\E(m)}$ for locally $\mathcal{E}(m)$-small maps in $\E(n)$. As in (\ref{eqn:51}), the  comonad 
 \[\f_n : \mathcal{E}(n) \to \mathcal{E}(n)\] 
 will induce, as in (\ref{eqn:52}), a comonad 
 \[\beta_m : F_n\pmb{\mho}_{\mathcal{E}(m)}  \to F_n\pmb{\mho}_{\mathcal{E}(m)} ,\]
 now of categories	in $\mathcal{E}^\f(n)$. 
	A universe $ \pmb{\mho}_{\mathcal{E}^\f(m)}$ for locally $\mathcal{E}^\f(m)$-small maps will then be constructed as the internal category of coalgebras $(F_n\pmb{\mho}_{\mathcal{E}(m)} )^{(\beta_m)}$. 

 \begin{lem}\label{lem:zwa23}
     Let $\F$ be a topos, $\mathrm{E}$ an impredicative universe in $\F$ in the sense of \cite{streicher2005universes}, and $\Box : \F \to \F$ a Cartesian comonad which takes $\mathrm{E}$-small maps to $\mathrm{E}$-small maps (equivalently which is display for  $\mathrm{E}$). Then, $\mathrm{E}^\Box$ is such an impredicative universe in $\F^\Box$.
 \end{lem}

 \begin{proof}
     This result appears as Corollary 5.2.1.3 of \citet{zwanziger2023natural}. To roughly trace the construction of a generic $\mathrm{E}^\f$-small map in $\F^\Box$, we proceed by choosing a generic $\mathrm{E}$-small map $p : \mho_\bullet \to \mho$ in $\F$. This is canonically extended to be the underlying action on objects of an internal functor $\pmb{p} : \pmb{\mho_\bullet} \to \pmb{\mho}$. (As usual, the internal category $\pmb{\mho}$ is the internal full subcategory associated to $p$.)

Because $p$ is generic, one can  form 
a pullback square 
\begin{equation*}
			\begin{tikzcd}
				\f \mho_\bullet \ar[dd, "\f p"']\ar[ddrr, phantom, "\lrcorner", very near start]\ar[rr, ""]&& \mho_\bullet \ar[dd, "p"]\\
				\\
				\f \mho  \ar[rr, ""'] &&\mho
			\end{tikzcd}
\end{equation*}
in $\F$. 
Transposing this in the adjunction $U \dashv F : \F \to \F^\Box$, one obtains a square 
\begin{equation*}
			\begin{tikzcd}
				F \mho_\bullet \ar[dd, "F p"']\ar[rr, "\beta_\bullet"]&& F \mho_\bullet \ar[dd, "F p"]\\
				\\
				F \mho  \ar[rr, "\beta"'] &&F\mho
			\end{tikzcd}
\end{equation*}
in $\F^\f$. These \[\beta : F \mho \to F \mho\] and \[\beta_\bullet : F \mho_\bullet \to F \mho_\bullet\] are canonically extended to be the underlying action on objects of internal comonads 
\[\pmb{\beta} : F \pmb{\mho} \to F \pmb{\mho}\] and \[\pmb{\beta_\bullet} : F \pmb{\mho_\bullet} \to F \pmb{\mho_\bullet},\]
respectively, such that the internal functor \[F \pmb{p}: F\pmb{\mho_\bullet} \to F \pmb{\mho}\] constitutes a strict morphism of comonads \[F \pmb{p} : \pmb{\beta_\bullet} \to \pmb{\beta}.\]
(In particular, the internal comonad \[\pmb{\beta} : F \pmb{\mho} \to F \pmb{\mho}\] is essentially an internalization of the indexed comonad 
\[\mathbb{B} : U^* \mathbb{E} \to U^* \mathbb{E}\]
induced by $\f$. One has in $\IndCat(\F^\Box)$ \begin{align*}
    [F \pmb{\mho}] & \simeq F_! [\pmb{\mho}] \\
                        & \simeq U^* [\pmb{\mho}] \\
                        & \simeq U^* \mathbb{E} ~,
\end{align*}
so $\mathbb{B}$ induces an indexed comonad on $[F \pmb{\mho}]$. Replacing its action on objects by $[\beta]$ strictifies this indexed comonad. The result then lies in the image of the externalization $2$-functor \[[-] : \Cat(\F^\Box) \monoto \Cat^{\F^{\Box}{\op}},\] yielding the internalization we call $\pmb{\beta}$.)

A generic $\mathrm{E}^\f$-small map is then obtained as the underlying action on objects of the internal functor \[  \Coalg(\pmb{\beta_\bullet}) \xrightarrow{\Coalg (F \pmb{p})} \Coalg(\pmb{\beta}),\]
constructed using the Cartesian structure of $\F^\Box$.

 \end{proof}

\begin{lem}\label{lem:topreprep}
Let $\E$ and $\F$ be small toposes, $\f_\E: \E \to \E$ and $\f_\F : \F \to \F$ be Cartesian comonads, and  \[ \f_\E \xhookrightarrow{i} \f_\F\]  a  dense logical inclusion of Cartesian comonads on toposes. Then, if the underlying dense logical inclusion of toposes \[\E \xhookrightarrow{i}\F\] admits a generic locally small map, so does the dense (by Lemma \ref{lem:topdensedense}) logical inclusion of toposes \[ \Coalg(\f_\E) \xhookrightarrow{\Coalg(i)} \Coalg(\f_\F).\] 
\end{lem}

\begin{proof} 
When we take the locally $\E$-small maps for display maps, the comonad $\f_\F : \F \to \F$ is display, as it is  Cartesian and takes locally $\E$-small maps to locally $\E$-small maps (Lemma \ref{lem:presls}). 
Thus, by Lemma \ref{lem:zwa23}, the set of display maps in $\Coalg(\f_\F)$ given by morphisms with an underlying locally $\E$-small map admits a generic display map. By the reasoning at Example \ref{ex:locsmall}, these display maps are precisely the locally $\Coalg(\f_\E)$-small maps. 
\end{proof}

\begin{defn}
When $\E$ 
and $\F$ are stratified toposes, we say that a morphism \[F : \E \to \F\] of $\cat{Topos}^\N$ is a \textbf{stratified Cartesian functor} 
if, for all $m \leq n$, the 
Cartesian functor \[ \E(n) \xrightarrow{F_n} \F(n)\]
constitutes a morphism of sites 
\[(\E(n), \E(m)) \xrightarrow{F_n} (\F(n),\F(m)).\]
\end{defn}

\begin{defn}
When $\E$ is a stratified topos, we say that a comonad $\f : \E \to \E$ in $\cat{Topos}^\N$ is a \textbf{stratified Cartesian comonad} if its underlying morphism is a stratified Cartesian functor.
\end{defn}

\begin{thm}\label{thm:stratcoalg}
Let $\E$ be a stratified topos and $\f : \E \to \E$ a stratified Cartesian comonad. Then, the evident functor $\E^\f : \N \to \Topos$ given on objects $n\in \N$ by $\E^\f(n) = \E(n)^{(\f_n)}$ is again a stratified topos. 
\end{thm}

\begin{proof}
By Lemma \ref{lem:topdensedense}, the inclusion functor $\E(m)^{(\f_m)} \leq \E(n)^{(\f_n)}$ is a dense logical inclusion, for all  $m \leq n$.
Moreover, by Lemma \ref{lem:topreprep}, the dense logical inclusion $\E(m)^{(\f_m)} < \E(n)^{(\f_n)}$ admits a generic locally 
small map, for all  $m < n$.   
\end{proof}

\begin{rem}
    Theorem \ref{thm:stratcoalg} provides most of the proof that the $2$-category of stratified toposes and stratified Cartesian functors admits the construction of coalgebras. However, we omit this statement due to it being of more technical interest.
\end{rem}

\section{Conclusion and Future Directions}\label{sec:conc}

We have presented a simple working definition of stratified topos. This is based on a definition of universe in a topos that strengthens \citet{streicher2005universes}'s notion of impredicative universe without excluding the main examples. By constructing the stratified topos of coalgebras, we have also shown that one of the key constructions of topos theory -- that of the topos of coalgebras -- can be refined to the present setting.

A loose end from the present work is the verification that realizability toposes yield examples of stratified
toposes in our sense, under some large cardinal assumptions. Alternative large cardinal axioms might also yield a realizability model of Streicher's axioms that contradicts our stratified topos axioms.

An avenue for future development would be to strengthen the present axioms. Some possibly desirable principles were already identified by work in algebraic set theory following \citet{joyal1995algebraic}. For example,  when $\E$ is a dense logical subtopos of $\F$, we might further ask that, in any pullback 

\begin{equation*}
\begin{tikzcd}
D\ar[dd, "o"' ]\ar[ddrr, phantom, "\lrcorner", very near start]\ar[rr, two heads, ""]&& E \ar[dd, "p"]\\
\\
A \ar[rr, two heads, "e"'] &&B
\end{tikzcd}
\end{equation*}
in $\F$ in which $e$ is an epimorphism,  the morphism $p$ is a locally small map whenever $o$ is.\footnote{Prospects for adding this axiom are discussed in \citet[][\S 4]{streicher2005universes}.} On another note, we might -- at least for some purposes -- require that universes in a stratified topos satisfy certain strict equations. For example, we might ask for an internal full and faithful functor \[\pmb{\mho}_m \xrightarrow{i} \pmb{\mho}_n\]
that preserves the internal topos structure of $\pmb{\mho}_m$ on the nose, for all  $m < n$. To accomplish this, we might pass to a natural model-theoretic approach, as in the natural display toposes of \citet{zwanziger2023natural}, or an approach imposing the realignment property, as in \citet{gratzer2022strict}.

Finally, there is much work to be done refining theorems of topos theory to the level of stratified toposes. We have seen an approach to constructing a surjection-embedding factorization for stratified geometric morphisms in our construction of the stratified topos of coalgebras. What about a hyperconnected-localic factorization? What is the right notion of open stratified geometric morphism?  Is a universe $\pmb{\mho}_n$ in a stratified topos $\E$ the `bi-terminal Grothendieck $\pmb{\mho}_n$-topos in $\E$,' just as the subobject classifier $\pmb{\Omega}$ in a topos $\E$ is the terminal locale in $\E$?

\section*{Acknowledgements}

I am grateful to Steve Awodey, Benno van den Berg, Peter Dybjer, Jonas Frey, Ansten Klev, Jaap van Oosten, Sam Speight, Thomas Streicher, and Andrew Swan for related feedback. A suggestion by Jonas led to my consideration of the right Kan extension functor $i_*$, which became central in this work. 

Work on this article was supported at various times by the project Logical Structure of Information Channels, No. 21-23610M of the Czech Science Foundation (carried out at the Institute of Philosophy of the Czech Academy of Sciences), the \textit{Lumina Quaeruntur} fellowship LQ300092101 from the Czech Academy of Sciences, and the (United States) Air Force Office of Scientific Research under award number FA9550-21-1-0009.

\appendix
\section{More General Coalgebra Constructions}\label{sec:appendix}

This section presents a few coalgebra constructions of general interest.

Recall that Lemma \ref{lem:topdensedense} gave conditions under which a dense subcategory 
\[\C \xhookrightarrow{i} \D\]
can be lifted to a dense subcategory of coalgebras 
\[ \C^{\f_\C} \inc \D^{\f_\D},\]
where $\f_\C$ and $\f_\D$ are  levelwise comonads. In particular, we required that $\C$ and $\D$ be toposes and $i$ be a logical inclusion. However, we used something much weaker, at most that $\C$, $\D$, and $i$ were Cartesian and that all split monos of $\D$ were locally $\C$-small. This last condition obtains, for example, when $\C$ contains a generic mono of $\D$. 

To lift such a generic mono to coalgebras, it is natural to let $\C$, $\D$, and $i$ be Heyting. One then makes use of the fact that the category of coalgebras for a Cartesian comonad on a Heyting category that admits a generic mono is again a Heyting category that admits a generic mono. This last is an apparent novelty, which we prove by applying the much more general results of \citet{zwanziger2023natural}.

Section \ref{sec:genmono} (resp.  Section \ref{sec:densesubcoal}) concerns lifting generic monos (resp. dense subcategories) to coalgebras.

\subsection{The Generic Mono of Coalgebras}\label{sec:genmono}

This section gives conditions under which a generic mono can be lifted to coalgebras. The construction is a natural generalization of that of the subobject classifier in the topos of coalgebras (Lawvere and Tierney, 1969 or 1970). Our approach is simply to apply \citet{zwanziger2023natural}.

\begin{proposition}\label{prop:heyt}
        When $\C$ is a Heyting category and $\Box : \C \to \C$ is a Cartesian comonad, $\C^\Box$ is again a Heyting category.
    \end{proposition}
    \begin{proof}
        See \citet[][the argument before Proposition 29]{warren2007coalgebras}.
    \end{proof}

    \begin{thm}\label{prop:heytgen}
        When $\C$ is a Heyting category that admits a generic monomorphism and $\Box : \C \to \C$ is a Cartesian comonad, $\C^\Box$ is again a Heyting category that admits a generic monomorphism.
    \end{thm}

    \begin{proof}
        In light of Proposition \ref{prop:heyt}, there remains to show that $\C^\Box$ admits a generic mono. This follows from Corollary 5.2.1.2 of \citet{zwanziger2023natural}, by taking the monos of $\C$ as the display maps.\footnote{The statement  in \citet{zwanziger2023natural} technically includes Cartesian closure, but this can simply be dropped.}
    \end{proof}

\begin{rem}
    As discussed in Section \ref{sec:stratcoalg}, the subobject  classifier in the topos of coalgebras is constructed as the internal poset of coalgebras (or, equivalently, the internal poset of fixed points) for the canonical induced internal comonad      $\beta : F\pmb{\Omega} \to F\pmb{\Omega}$, where $\pmb{\Omega}$ denotes the subobject classifier of the underlying topos, viewed as an internal poset, and $F$ is the cofree functor's action on internal posets. This can familiarly be calculated as the equalizer  of  $\id_{F\pmb{\Omega}}$ and $\beta$.

    In the construction used in Theorem \ref{prop:heytgen}, we  similarly construct the generic mono (in its guise as a weak subobject classifier) as the internal preorder of coalgebras for the canonical induced internal comonad      $\beta : F\pmb{\Omega} \to F\pmb{\Omega}$, where now $\pmb{\Omega}$ denotes the weak subobject classifier of $\C$, viewed as an internal preorder, and $F$ is the cofree functor's action on internal preorders. 

    However, it is not appropriate to calculate the latter object of coalgebras as the equalizer  of  $\id_{F\pmb{\Omega}}$ and $\beta$, since coalgebras for a comonad on a preorder need not be fixed points.   Rather, we may take the inserter  of  $\id_{F\pmb{\Omega}_{\C}}$ and $\beta$. 
\end{rem}

We will also make use of the generalization from Heyting categories to inclusions thereof. 

In a context where $\C$ and $\D$ are each a Heyting category that admits a generic mono, we will say that a Cartesian full subcategory inclusion $i: \C \inc \D$ is a \textbf{logical inclusion} if it is Heyting and preserves the generic mono.

And where $\C$ and $\D$ are each a Heyting category that admits a generic mono and $\f_\C: \C \to \C$ and $\f_\D : \D \to \D$ are Cartesian comonads,  we will say that a Cartesian full inclusion of Cartesian comonads  \[\f_\C \xhookrightarrow{i} \f_\D\]  
 is a \textbf{logical inclusion}  of Cartesian comonads
if the underlying Cartesian full subcategory inclusion \[\C \xhookrightarrow{i}  \D\] is a logical inclusion.

\begin{lem}\label{lem:comonloglogh}
When $\C$ and $\D$ are each a Heyting category that admits a generic mono, $\f_\C: \C \to \C$ and $\f_\D : \D \to \D$ are Cartesian comonads, and  \[\f_\C \xhookrightarrow{i} \f_\D\] is a  logical inclusion of Cartesian comonads, the Cartesian full subcategory inclusion \[ \Coalg(\f_\C) \xhookrightarrow{\Coalg(i)} \Coalg(\f_\D)\] is a  logical inclusion.
\end{lem}

\begin{proof}
The underlying logical inclusion $i : \C \inc \D$ preserves and commutes with everything in sight, including the construction of Heyting category structure and generic monos  in coalgebras. The inclusion $\Coalg(i) : \Coalg(\f_\C) \to \Coalg(\f_\D)$ is thus logical.
\end{proof}

\subsection{The Dense Subcategory of Coalgebras}\label{sec:densesubcoal}

This section gives conditions under which a dense subcategory can be lifted to coalgebras. 

When $\C$ and $\D$ are Cartesian categories  and $\f_\C: \C \to \C$ and $\f_\D : \D \to \D$ are Cartesian comonads, we say that a Cartesian full inclusion of Cartesian comonads \[\f_\C \inc \f_\D\] is  \textbf{dense} if its underlying  Cartesian full subcategory inclusion \[\C \inc \D\] is a dense functor and the underlying Cartesian functor 
\[ \D \xrightarrow{\f_\D} \D\]
constitutes a morphism of sites
\[ (\D,\C) \xrightarrow{\f_\D} (\D,\C).\]

\begin{thm}\label{lem:densedense}
Let $\C$ and $\D$ be Cartesian categories, $\f_\C: \C \to \C$ and $\f_\D : \D \to \D$ be Cartesian comonads, and  \[ \f_\C \inc \f_\D\] be a  dense Cartesian full inclusion of Cartesian comonads. If all split monomorphisms of $\D$ are locally $\C$-small, then all split monomorphisms of $\D^{\f_\D}$ are locally $\C^{\f_\C}$-small and the Cartesian full subcategory inclusion \[ \C^{\f_\C} \inc \D^{\f_\D}\] is dense.
\end{thm}

\begin{proof}
    In the argument of Example \ref{ex:locsmall}, we used only that all split monomorphisms of $\D$ are locally $\C$-small to establish that the forgetful functor \[U_\D : \D^{\f_\D}  \to \D\]
constitutes a morphism of sites
\[U_\D : (\D^{\f_\D},\C^{\f_\C})  \to (\D,\C),\]
and, therefore, that the morphisms of $\D^{\f_\D}$ with underlying locally $\C$-small map coincide with the locally $\C^{\f_\C}$-small maps. 

The argument of Lemma \ref{lem:topdensedense} then establishes that the Cartesian full subcategory inclusion $\C^{\f_\C} \inc \D^{\f_\D}$ is dense.

Any split monomorphism of $\D^{\f_\D}$ is locally $\C^{\f_\C}$-small, as its underlying split monomorphism is locally $\C$-small by assumption.  
\end{proof}

\begin{thm}

Let $\C$ and $\D$ each be a Heyting category that admits a generic mono, $\f_\C: \C \to \C$ and $\f_\D : \D \to \D$ be Cartesian comonads, and  \[\f_\C \xhookrightarrow{i} \f_\D\] be a dense  Cartesian full inclusion of Cartesian comonads, which is, moreover,  a logical inclusion. Then, the categories of coalgebras $\C^{\f_\C}$ and $\D^{\f_\D}$ are again each a Heyting category that admits a generic mono, and  the Cartesian full subcategory inclusion \[ \C^{\f_\C} \inc \D^{\f_\D}\] is  dense and a logical inclusion.

			\end{thm}

   \begin{proof}
        The argument of Proposition \ref{lem:logsubmon}, requires only the assumption of a generic mono in $\C$ that is preserved by the inclusion to $\D$ to establish that all monomorphisms of $\D$ are locally $\C$-small. This includes the split monomorphisms in particular.

We can thus apply Lemma \ref{lem:comonloglogh} and Theorem \ref{lem:densedense}.
   \end{proof}

\bibliographystyle{plainnat}
\bibliography{bib} 

\end{document}